\documentclass[11pt,twoside,a4paper]{amsart}
\usepackage{color}
\usepackage[colorlinks=true,urlcolor=blue,linkcolor=blue,citecolor=black]{hyperref}
\usepackage{amscd}
\usepackage{amsmath}
\usepackage{amssymb}
\usepackage{amsthm}
\usepackage{latexsym}
\usepackage{stmaryrd}
\usepackage{float}
\usepackage{shuffle}
\usepackage{tikz,tikz-cd}

 \usepackage{hyperref}
\hypersetup{
    colorlinks=true,
    linkcolor=blue,
    filecolor=magenta,      
    citecolor=green,
    urlcolor=cyan,
    pdftitle={Pre Lie Algebras, Forest Formula and Cumulants},
    pdfpagemode=FullScreen,
    }

\usepackage[T1]{fontenc}   
\usepackage{comment}
\usepackage[english]{babel}

\theoremstyle{definition}
\newtheorem{defi}{Definition}[section]
\newtheorem{rem}[defi]{Remark}
\newtheorem{ex}[defi]{Example}
\newtheorem{nota}[defi]{Notation}
\theoremstyle{theorem}
\newtheorem{lemma}[defi]{Lemma}

\newtheorem{theo}[defi]{Theorem}
\newtheorem{prop}[defi]{Proposition}

\newcommand{\frontstick}{\,\raisebox{-1pt}{\begin{tikzpicture}
\draw [line width=1pt,] (0,0)--(0,0.25);
\end{tikzpicture}}\kern+2pt}

\newcommand{\Id}{\operatorname{Id}}

\newcommand{\sol}{\operatorname{sol}_1}

\newcommand{\A}{\mathcal{A}}

\newcommand{\NC}{\operatorname{NC}}

\newcommand{\T}{\mathcal{T}}

\newcommand{\Lfr}{L_{\mathrm{free}}}

\newenvironment{arb}{\begin{tikzpicture}[baseline,scale=0.5,level distance=7mm,level 1/.style={sibling distance=10mm},level 2/.style={sibling distance=5mm},level 3/.style={sibling distance=3mm},grow=down, font=\scriptsize]
\tikzstyle{ve}=[draw,circle,inner sep=1pt,fill] 
\tikzstyle{vv}=[draw,circle,inner sep=1pt] 
\tikzstyle{vee}=[minimum size=0pt ,inner sep=0pt]}{\end{tikzpicture}}
\newenvironment{arbb}{\begin{tikzpicture}[baseline,scale=0.5,level distance=7mm,level 1/.style={sibling distance=25mm},level 2/.style={sibling distance=5mm},level 3/.style={sibling distance=3mm},grow=down, font=\scriptsize]
\tikzstyle{ve}=[draw,circle,inner sep=1pt,fill] 
\tikzstyle{vv}=[draw,circle,inner sep=1pt] 
\tikzstyle{vee}=[minimum size=0pt ,inner sep=0pt]}{\end{tikzpicture}}

\newcommand{\rd}[1]{\node[ve,label=above:$#1$] {}}

\newcommand{\vl}[1]{node[ve,label=left:$#1$] {}}
\newcommand{\vr}[1]{node[ve,label=right:$#1$] {}}
\newcommand{\vb}[1]{node[ve,label=below:$#1$] {}}

\def\corollaa{\begin{arb}
\rd{(i;i_0)}
child{\vb{i_1}} child[dotted] child[dotted] child{\vb{i_k}};
\end{arb}}

\def\corollab{\begin{arb}
\rd{(i;i_0)}
child{\vb{i_1}}  child{\vb{i_2}};
\end{arb}}

\def\corollac{\begin{arb}
\rd{(i;i_0)}
child{\vl{(i_1 ; i_{1,0})} child{\vb{i_{1,1}}} child[dotted] child[dotted] child{\vb{i_{1,k}}}} child{\vr{i_2}};
\end{arb}}

\def\corollad{\begin{arb}
\rd{(i;i_0)}
child{\vl{i_1}} child{\vr{(i_2 ; i_{2,0})}  child{\vb{i_{2,1}}} child[dotted] child[dotted] child{\vb{i_{2,l}}}} ;
\end{arb}}

\def\corollae{\begin{arbb}
\rd{(i;i_0)}
child{\vl{(i_1 ; i_{1,0})} child{\vb{i_{1,1}}} child[dotted] child[dotted] child{\vb{i_{1,k}}}}child{\vr{(i_2 ; i_{2,0})}  child{\vb{i_{2,1}}} child[dotted] child[dotted] child{\vb{i_{2,l}}}} ;
\end{arbb}}

\usepackage{forest}
\forestset{
  decor/.style = {
    label/.expanded = {[inner sep = 0.2ex, font=\unexpanded{\tiny}]right:{$#1$}}
  },
  root/.style = {minimum size = 0.1ex},
  decorated/.style = {
    for tree = {
      circle, fill, inner sep = 0.3ex, minimum size = 1.ex,
      grow' = south, l = 0, l sep = 1.2ex, s sep = 0.7em,
      fit = tight, parent anchor = center, child anchor = center,
      delay = {decor/.option = content, content =}
    }
  },
  default preamble = {decorated, root},
  begin draw/.code={\begin{tikzpicture}[baseline={([yshift=-0.5ex]current bounding box.center)}]},
}

\newcommand{\NCirr}{\operatorname{NC}_{\mathrm{irr}}}

\definecolor{red}{rgb}{1.,0.,0.}
\definecolor{green}{rgb}{0.,1.,0.}
\definecolor{blue}{rgb}{0.,0.,1.}
\definecolor{orange}{rgb}{1.,0.8431372549019608,0.}

\title[A Forest formula for pre-Lie exponentials, Magnus'operator and cumulants]{A forest formula for pre-Lie exponentials, Magnus' operator and cumulant-cumulant relations}

\author{Adrián Celestino and Frédéric Patras}
\address{Department of Mathematical Sciences, Norwegian University of Science and Technology (NTNU), 7491 Trondheim, Norway.}
\email{adrian.celestino@ntnu.no}
\address{Universit\'e C\^ote d'Azur, Laboratoire J.-A.~Dieudonn\'e, UMR 7351, CNRS, Parc Valrose, 06108 Nice Cedex 02, France.}
\email{frederic.patras@unice.fr}

\begin{document}



\tikzset{
solid node/.style={circle,draw,inner sep=1.5,fill=black, minimum size=1pt},
hollow node/.style={circle,draw,inner sep=1.5,fill=white,  minimum size=1pt}
}

\maketitle
\begin{abstract}
    Forest formulas that generalize Zimmermann's forest formula in quantum field theory have been obtained for the computation of the antipode in the dual of enveloping algebras of pre-Lie algebras. In this work, largely motivated by Murua's analysis of the Baker-Campbell-Hausdorff formula, we show that the same ideas and techniques generalize and provide effective tools to handle computations in these algebras, which are of utmost importance in numerical analysis and related areas. We illustrate our results by studying the action of the pre-Lie exponential and the Magnus operator in the free pre-Lie algebra and in a pre-Lie algebra of words originating in free probability. The latter example provides combinatorial formulas relating the different brands of cumulants in non-commutative probability.
\end{abstract}

\section{Introduction}
In the theory of elementary particles (quantum field theory, QFT), one often needs to remove infinities to get sound quantities. This process is called renormalization and involves computing a {\it counterterm}. Mathematically, a counterterm is the convolution inverse of a character on a certain Hopf algebra $H$ of Feynman diagrams  (an algebra map from $H$ to the scalars). Renormalization is actually a general process that has been used recently in a wide variety of contexts such as number theory and stochastic processes.

Computing the convolution inverse of a character on a Hopf algebra boils down to computing the antipode of the Hopf algebra (the convolution inverse of the identity). In QFT, this is achieved by the so-called  Zimmermann's forest formula. With respect to other approaches, Zimmermann's forest formula has the nice property to reduce
drastically the number of terms in the expression of the antipode (it is optimal in that sense \cite{FG05}). 
The Hopf algebras involved in renormalization have the specific feature of being the dual to enveloping algebras of pre-Lie algebras. It was shown in \cite{mp} that Zimmermann's forest formula generalizes and is best understood in this broader framework. We refer to the book \cite{CP21} for general results on Hopf and pre-Lie algebras and to its Chap. 10 for an introduction to renormalization.

The present article was motivated by our recent work on cumulant-cumulant relations in non-commutative probability \cite{cefpp21} and the work of Murua \cite{Mu06}. Murua's work is a key contribution to one of the most classical problems: the Baker-Campbell-Hausdorff (BCH) problem, that is the computation of the logarithm of the solution of a non-autonomous matrix-valued linear differential equation. In group theory, this is the problem of giving explicit formulas relating the group law of a Lie group with the Lie product law of its Lie algebra. On another hand, cumulant-cumulant formulas have also be proven to be of group/Lie theoretical nature \cite{EFP19}.

In our \cite{cefpp21}, the same combinatorics of trees and the same combinatorial coefficients as in \cite{Mu06} appeared. The reason for this appearance was easy to guess: in non-commutative probability, a pre-Lie structure underlies our formulas, whereas trees encode the free pre-Lie algebra on one generator.
This prompted us to look for a more systematic understanding of the formulas in these works using pre-Lie algebras techniques. 
This goal is achieved here. In the process, we develop a new strategy for the computation of the Magnus operator, which is based on the particular form taken by the canonical projection from the enveloping algebra of a pre-Lie algebra onto the underlying pre-Lie algebra \cite{CHP13}.

Concretely, we first show that the technique leading to a general pre-Lie forest formula for the antipode in the dual of the enveloping algebra of a pre-Lie algebra can be adapted to get closed formulas for the computation of iterated pre-Lie products. This is a general and effective process. Effective means here that, once the structure constants of a pre-Lie algebra (for example of vector fields) in a basis are known, it can be applied to get explicit formulas for iterated pre-Lie products. In the second part of the article, we first apply these results to two fundamental operations in the pre-Lie theory: the pre-Lie exponential, or Agrachev-Gamkrelidze operator, and its inverse, the Magnus operator (see \cite{AG} and \cite[Chap. 6]{CP21}). In the context of differential equations, the Magnus formula, from which the name ``Magnus operator'' was derived, solves the BCH problem: this is the connection to \cite{Mu06}.
Ultimately, we apply these results to free probability and cumulant-cumulant formulas. In the process, we correct one of the coefficients appearing in the pre-Lie Zimmermann's forest formula of \cite{mp}, where a symmetry factor was missing.
\subsection*{Organization of the paper} In addition to the preceding introductory section, this paper is organized as follows. Section \ref{sec:preliebrace} contains the definitions of the algebraic structures used in this manuscript: pre-Lie and symmetric brace algebras, as well as the connections between them. In Section \ref{sec:preliewords}, we expand the example of a pre-Lie algebra of words originating in noncommutative probability, and obtain formulas for the symmetric brace product and the dual coproduct associated to the pre-Lie algebra structure. In Section \ref{sec:iteratedprelie}, we state and prove some technical lemmas related to the action on iterated pre-Lie and iterated brace products of elements in the dual space. Section \ref{sec:forestformulas} is devoted to obtain forest formulas for different types of iterated coproducts. In particular, we correct the forest formula appearing in \cite{mp} by adding a symmetry coefficient associated to every decorated tree. The subsequent sections correspond to applications of the machinery developed in this work. Section \ref{sec:preliexexpfree} details the computation of the pre-Lie exponential of the generator of the free pre-Lie algebra of rooted trees by using the forest formulas for iterated coproducts. In Sections \ref{sec:MuruaMagnusI} and \ref{sec:MuruaMagnusII}, we show how to obtain the action of the Magnus operator on the generator of the free pre-Lie algebra following two different approaches relying purely on pre-Lie algebra techniques. The inclusion of two different proofs of the same formulas was motivated by the importance of the subject, and, for the second, by the goal of making explicit connections with the existing literature. Finally, in Section \ref{sec:cumulantformulas}, we present an application of our methods in the context of combinatorial relations, in terms of non-crossing partitions, between monotone, free, and Boolean cumulants in non-commutative probability.

\section{Pre-Lie and Symmetric Brace algebras}
\label{sec:preliebrace}

In this section, we recall first generalities on coalgebras and pre-Lie algebras and their enveloping algebras. The reader can find a more detailed presentation in \cite[Chap. 6]{CP21}.

Let $V$ be a vector space over a ground field $\kappa$ of characteristic $0$. Recall that the cofree cocommutative connected graded coalgebra over $V$ identifies with $\kappa[V]$, the vector space of polynomials over $V$.
We write $\kappa[V]_n$ for the linear span of degree $n$ monomials, so that, in tensor notation, $\kappa[V]_n=(V^{\otimes n})_{S_n}$, where $S_n$ stands for the $n$-th symmetric group.
Under this identification, the coproduct can be written as
$$\Delta(l_1\cdots l_n) = \sum_{I\subseteq [n]} l_I\otimes l_{[n]\backslash I},$$
where $l_1,\ldots,l_n \in V$ and, given of subset $S=\{s_1,\ldots,s_k\}$ of $[n]:=\{1,\ldots,n\}$,  $l_S:=l_{s_1}\cdots l_{s_k}$. The graduation is the usual one: it is induced by the degrees of monomials. The counit $\varepsilon$ is the projection on the scalar component: the map sending a polynomial to its constant term.
The product of polynomials is written $\cdot$, and it makes the triple $(\kappa[V],\cdot,\Delta)$ a Hopf algebra, that is, the coproduct $\Delta$ and the counit $\varepsilon$ are maps of unital algebras from $\kappa[V]$ to $\kappa[V]\otimes \kappa[V]$. The antipode is the algebra automorphism of $\kappa[V]$ induced by the linear automorphism of $V$, $v\mapsto -v$.

The duality pairing $V\otimes V^*\to k$ denoted $\langle v|w^* \rangle, \ v\in V,\ w^\ast\in V^\ast$, with a self-explaining notation, extends to a pairing $\langle \,\cdot\, |\,\cdot\, \rangle :\kappa[V]\otimes \kappa[V^*]\to \kappa$, defined by requiring $\langle P|Q\rangle =0$ if $P$ and $Q$ are monomials in $\kappa[V]$ resp.  $\kappa[V^\ast]$ of distinct degrees and:

\begin{align} 
\label{eq:extension}
         \langle v_1\cdots v_n|w_1\cdots w_n \rangle &=\sum_{\sigma\in S_n}\prod\limits_{i=1}^n \langle v_{\sigma(i)}|w_i\rangle,
\end{align}
for any $n\geq1$, $v_1,\ldots,v_n\in V,$ and $w_1,\ldots,w_n\in V^*$. We also recall the notion of graded dual: if $W=\bigoplus_{n\geq1}W_n$ is a graded vector space, then its graded dual is $W^\ast:=\bigoplus_{n\geq1 }W_n^\ast$ (see e.g. \cite{CP21} for details and related notions relevant to the present article --such as graded filtered or complete vectors spaces and Hopf algebras). 
When $V$ is finite dimensional, the above pairing allows, in particular, to identify $\kappa[V^*]$ with the graded dual of $\kappa[V]$. 

Similar observations hold when $V$ is a locally finite connected graded vector space: that is, a graded vector space such that each graded component is finite dimensional and $V_0=0$. The latter condition corresponds to {connectedness} for graded vector spaces. To take into account units, a graded algebra $A$ is said to be \textit{connected} if it holds instead that $A_0=\kappa$. In that case, $\kappa[V]$ is equipped with two graduations: the ``polynomial'' one obtained from the degrees of monomials, and the ``full'' graduation induced by the graduation of $V$, so that a product $v_1\cdots v_n$ of elements $v_i\in V_{p_i}$ is of degree $p_1+\cdots + p_n$. For a locally finite connected graded vector space $V$, the graded dual of $\kappa[V]$ is defined using the second (the full) graduation.

For later use, given an arbitrary coassociative coalgebra $C$ equipped with a coproduct $\Delta$, we introduce the Sweedler notation, $\Delta(P)=:P^{(1)}\otimes P^{(2)}$, the iterated coproduct, $\Delta^{[1]}:=\Id$, $\Delta^{[n+1]}:=(\Delta\otimes \Id^{\otimes n-1})\circ \Delta^{[n]}$, from $C$ to $C^{\otimes n+1}$, $n\geq 1$, and the generalized Sweedler notation,
$\Delta^{[n]}(P)=:P^{(1)}\otimes\dots\otimes P^{(n)}$.

\begin{defi}
\label{def:pre-Lie}
A \textit{pre-Lie algebra} is a vector space $L$ equipped with a linear map $\lhd: L\otimes L\to L$, $v\otimes w\mapsto v\lhd w$ such that
$$(u\lhd v)\lhd w - u \lhd (v\lhd w) = (u\lhd w)\lhd v - u \lhd (w\lhd v).$$
\end{defi}

Given a pre-Lie algebra $L$, consider the coalgebra $\kappa[L]$. It can be equipped with another associative product $*$, making $(\kappa[L],*,\Delta)$ a Hopf algebra and the enveloping algebra of $L$. This product $*$ identifies with $\cdot$ if and only if $\lhd$ is the null product.
More precisely, we have
\begin{equation}
    \label{eq:product}
     (a_1\ldots a_\ell)*(b_1\cdots b_m) = \sum_f B_0(a_1\{B_1\})\cdots (a_\ell \{B_\ell\}),
\end{equation}
where the sum is over all the maps $f:\{1,\ldots,m\}\to \{0,\ldots \ell\}$, 
$B_i:=b_{f^{-1}(\{i\})}$,
and the meaning of the braces is defined as follows.

\begin{defi}
A \textit{symmetric brace algebra} is a vector space $V$ equipped with a linear map $V\otimes \kappa[V] \to V$, $v\otimes P\mapsto v\{P\}$ such that $v\{1\} = v$ and 
$$ (v\{v_1,\ldots,v_n\})\{P\} = v\{v_1\{P^{(1)}\},v_2\{P^{(2)}\},\ldots,v_n\{ P^{(n)}\},P^{(n+1)}\},$$
where $\Delta^{[n+1]}(P) = P^{(1)}\otimes \cdots \otimes P^{(n+1)}.$
\end{defi}

\begin{theo}[Oudom-Guin \cite{OudomGuin} {\cite[Thm. 6.2.1]{CP21}}]
\label{thm:OudomGuin}
A pre-Lie algebra $(L,\lhd)$ is equipped with the structure of a symmetric brace algebra by the formulas
\begin{equation}
\label{eq:Oudom}
\begin{split}
v\{w\} &= v\lhd w,  \\  v\{w_1,\ldots,w_n\} &= (v\{w_1,\ldots,w_{n-1}\})\{w_n\} - \sum_{i=1}^{n-1} v\{w_1,\ldots,w_i\{w_n\},\ldots,w_{n-1}\},
\end{split}
\end{equation}
for all $v,w_1,\ldots,w_n\in L$. 
\end{theo}

Observe for later use that, when $L$ is finite dimensional or connected graded and locally finite dimensional, the brace product $L\otimes \kappa[L]_n\to L$ dualizes to a linear map ${\delta}_n:L^*\to L^*\otimes \kappa[L^*]_n$. Take care that in the graded and locally finite dimensional case, the index $n$ stands here for the {\it polynomial degree}. That is, a product $\lambda_1\cdots \lambda_n$ of elements of $L^\ast$ is of polynomial degree $n$, whatever the degrees of the $\lambda_i$ in $L^\ast$ -- whereas its full degree is, if the $\lambda_i$ are homogeneous, the sum of the degrees of the $\lambda_i$. Notice that, for degree reasons, ${\delta}_n=0$ on $L_i^\ast:=(L_i)^\ast$ for $n\geq i$. Hence, if we define $\overline{\delta} = \sum_{n\geq1} \delta_n$ in $L^*$ and $\delta:\kappa[L^*] \to \kappa[L^*]\otimes \kappa[L^*]$ to be the algebra morphism given by $\delta(w) = \overline{\delta}(w) + 1\otimes w + w\otimes 1$ for any $w\in L^*$, we have the following result.

\begin{theo}
Let $L$ be a connected graded and locally finite dimensional pre-Lie algebra. Then $(\kappa[L^*],\cdot, \delta)$ defines the Hopf algebra structure dual to the Hopf algebra $(\kappa[L], *, \Delta)$.
\end{theo}

\section{A Pre-Lie Algebra of Words}
\label{sec:preliewords}

We now introduce an important example of pre-Lie algebras for our purposes. Its introduction is motivated by the study of cumulant-cumulant relations in non-commutative probability: see our \cite{cefpp21} and the forthcoming developments in this article. Henceforth, let $X$ be a finite alphabet and denote by $X^*$ the set of \textit{non-empty} words over $X$. In addition, consider the vector space $L = \kappa X^*$ given by the linear span of the set of {non-empty} words over $X$. It is a graded vector space, $L=\bigoplus_{n\geq1}L_n$, where $L_n$ denotes the linear span of words of length $n$ for $n\geq 1$ and $L_0:=0$. Finally, given two words $\alpha_1,\alpha_2\in X^*$, we denote by $\alpha_1\alpha_2$ to the concatenation of $\alpha_2$ to the right of $\alpha_1$. Now, we define the linear map $\lhd: L\otimes L \to L$ by
\begin{equation}
\label{eq:pre-LieWords}
\alpha\lhd \gamma = \sum_{\substack{\alpha_1\alpha_2 = \alpha\\\alpha_1,\alpha_2\neq \emptyset}} \alpha_1\gamma\alpha_2,\qquad \forall\, \alpha,\gamma\in X^*.
\end{equation}

This map is graded: $\alpha\lhd\gamma\in L_{n+m}$ if $\alpha\in L_n$ and $\gamma\in L_m$.

\begin{prop}
Let $L = \kappa X^*$ be the vector space given by the linear span of $X^*$. Then $(L,\lhd)$ is a graded pre-Lie algebra, with $\lhd$ given in \eqref{eq:pre-LieWords}.
\end{prop}
\begin{proof}
We will prove that the pre-Lie identity in Definition \ref{def:pre-Lie} holds. With this purpose, take words $\alpha,\gamma,\xi\in X^*$. Observe that
\begin{eqnarray*}
(\alpha\lhd \gamma) \lhd \xi &=& \left( \sum_{\substack{\alpha_1\alpha_2 = \alpha\\\alpha_1,\alpha_2\neq \emptyset}} \alpha_1\gamma\alpha_2 \right)\lhd \xi
\\ &=& \sum_{\substack{\alpha_1\alpha_2 = \alpha\\\alpha_1,\alpha_2\neq \emptyset}} \left( \sum_{\substack{\alpha_{1,1}\alpha_{1,2}= \alpha_1\\ \alpha_{1,1}\neq\emptyset}} \alpha_{1,1}\xi\alpha_{1,2}\gamma\alpha_{2} + \sum_{\substack{\gamma_1\gamma_2=\gamma\\ \gamma_1,\gamma_2\neq\emptyset}} \alpha_1\gamma_1\xi\gamma_2\alpha_2
+\sum_{\substack{\alpha_{2,1}\alpha_{2,2}= \alpha_2\\ \alpha_{2,2}\neq\emptyset}} \alpha_1\gamma\alpha_{2,1}\xi\alpha_{2,2} \right).
\end{eqnarray*}
On the other hand
\begin{eqnarray*}
\alpha \lhd( \gamma\lhd \xi) = \alpha \lhd \left( \sum_{\substack{\gamma_1\gamma_2 = \gamma\\ \gamma_1,\gamma_2\neq\emptyset}} \gamma_1\xi\gamma_2\right) = \sum_{\substack{\alpha_1\alpha_2 = \alpha \\ \alpha_1,\alpha_2\neq\emptyset}}\sum_{\substack{\gamma_1\gamma_2 = \gamma\\ \gamma_1,\gamma_2\neq\emptyset}} \alpha_1\gamma_1\xi\gamma_2\alpha_2
\end{eqnarray*}
Therefore
\begin{equation*}
(\alpha\lhd \gamma) \lhd \xi - \alpha \lhd( \gamma\lhd \xi) =  \sum_{\substack{\alpha_1\alpha_2\alpha_3 = \alpha\\\alpha_1,\alpha_3\neq \emptyset}} \left(  \alpha_{1}\xi\alpha_{2}\gamma\alpha_{3} + \alpha_1\gamma\alpha_{2}\xi\alpha_{3} \right).
\end{equation*}
As the expression in the right hand side is symmetric in $\gamma$ and $\xi$, we  obtain that it also computes
$$(\alpha\lhd \xi) \lhd \gamma - \alpha \lhd( \xi\lhd \gamma). $$ from which we conclude that the pre-Lie identity holds.
\end{proof}

By Theorem \ref{thm:OudomGuin}, $\lhd$ gives rise a symmetric brace operation on $L$ via Equations \eqref{eq:Oudom}. We are now interested in finding a closed formula for this map.

\begin{prop}
\label{prop:braceproduct}
Let $L = \kappa X^*$ be the pre-Lie algebra of words over $X$. The symmetric brace map given by Theorem \ref{thm:OudomGuin} is given by
\begin{equation}
\label{eq:braceL}
\alpha\{\gamma_1,\ldots,\gamma_n\} = \sum_{\sigma\in S_n} \sum_{\substack{\alpha_1\cdots \alpha_{n+1} = \alpha\\ \alpha_1,\alpha_{n+1}\neq \emptyset}} \alpha_1 \gamma_{\sigma(1)} \alpha_2 \gamma_{\sigma(2)} \cdots \alpha_n \gamma_{\sigma(n)} \alpha_{n+1},
\end{equation}
for any $\alpha,\gamma_1,\ldots,\gamma_n\in X^*$, $n\geq1$.
\end{prop} 
\begin{proof}
We will prove the statement by induction on $n$. The base case $n=1$ is straightforward from the definition of the pre-Lie and symmetric product. Now assume that the formula \eqref{eq:braceL} holds for a positive integer $n=k\geq1$. We aim to prove that the formula \eqref{eq:braceL} also holds for $n=k+1$.  In order to do this, we take $\alpha,\gamma_1,\ldots,\gamma_{k+1}\in X^*$. Thus we have that
\begin{eqnarray*}
\alpha\{\gamma_1,\ldots,\gamma_{k+1}\} &=& \left(\alpha\{\gamma_1,\ldots,\gamma_{k}\} \right)\{\gamma_{k+1}\} - \sum_{i=1}^{k} \alpha\{\gamma_1,\ldots,\gamma_i\{\gamma_{k+1}\},\ldots,\gamma_k\}
\end{eqnarray*}
By the inductive hypothesis, the first term of the right-hand side of the previous equation is equal to
\begin{eqnarray*}\left(\alpha\{\gamma_1,\ldots,\gamma_{k}\} \right)\{\gamma_{k+1}\} = \left( \sum_{\sigma\in S_k} \sum_{\substack{\alpha_1\cdots \alpha_{n+1} = \alpha\\ \alpha_1,\alpha_{k+1}\neq \emptyset}} \alpha_1 \gamma_{\sigma(1)} \alpha_2 \gamma_{\sigma(2)} \cdots \alpha_k \gamma_{\sigma(k)} \alpha_{k+1}\right)\{\gamma_{k+1}\} 
\end{eqnarray*}
According to the definition of the pre-Lie product, $\alpha_1 \gamma_{\sigma(1)}  \cdots \alpha_k \gamma_{\sigma(k)} \alpha_{k+1}\lhd \gamma_{k+1}$ is a sum indexed by the subwords generated by inserting $\gamma_{k+1}$ in $\alpha_1 \gamma_{\sigma(1)}  \cdots \alpha_k \gamma_{\sigma(k)}$. This process induces a partition of every word $\alpha_i$ and $\gamma_j$ into two subwords, with the only one restriction that the left subword of $\alpha_1$ and the right subword of $\alpha_{k+1}$ are non-empty. In other words, we have
\begin{eqnarray*}
\nonumber \left(\alpha\{\gamma_1,\ldots,\gamma_{k}\} \right)\{\gamma_{k+1}\} &=& \sum_{\sigma\in S_k}  \sum_{\substack{\alpha_1\cdots \alpha_{k+1} = \alpha\\ \alpha_1,\alpha_{k+1}\neq \emptyset}} 
\left( \sum_{\substack{\alpha_{1,1}\alpha_{1,2}=\alpha_1\\ \alpha_{1,1}\neq\emptyset }}  \alpha_{1,1}\gamma_{k+1}\alpha_{1,2}\gamma_{\sigma(1)}\cdots \gamma_{\sigma(k)} \alpha_{k+1}\right. \\ 
& &+\left. \sum_{i=2}^{k} \sum_{\substack{\alpha_{i,1}\alpha_{i,2}=\alpha_i\\\phantom{}}} \alpha_1\gamma_{\sigma(1)}\cdots \alpha_{i,1}\gamma_{k+1}\alpha_{i,2}\cdots \gamma_{\sigma(k)} \alpha_{k+1}\right. \\ 
& &+\left. \sum_{\substack{\alpha_{k+1,1}\alpha_{k+1,2}=\alpha_{k+1}\\
\alpha_{k+1,2}\neq\emptyset }}
\alpha_{1}\gamma_{\sigma(1)}\cdots \alpha_{k+1,1}\gamma_{k+1}\alpha_{k+1,2}\right. \\ 
& &  +\left. \sum_{i=1}^k\sum_{\substack{\gamma_{i,1}\gamma_{i,2}= \gamma_{\sigma(i)}\\ \gamma_{i,1},\gamma_{i,2}\neq\emptyset}}  \alpha_1\cdots \alpha_i\gamma_{i,1}\gamma_{k+1}\gamma_{i,2}\alpha_{i+1}\cdots \alpha_{k+1} \right)
\end{eqnarray*}
On the other hand,
\begin{eqnarray*}
\sum_{i=1}^{k} \alpha\{\gamma_1,\ldots,\gamma_i\{\gamma_{k+1}\},\ldots,\gamma_k\} &=& \sum_{i=1}^k \sum_{\substack{\gamma_{i,1}\gamma_{i,2} = \gamma_i\\ \gamma_{i,1},\gamma_{i,2}\neq\emptyset}} \alpha\{ \gamma_1,\ldots,\gamma_{i,1}\gamma_{k+1}\gamma_{i,2},\ldots,\gamma_k\}
\\ &=&  \sum_{i=1}^k \sum_{\substack{\gamma_{i,1}\gamma_{i,2} = \gamma_i\\ \gamma_{i,1},\gamma_{i,2}\neq\emptyset}} \sum_{\sigma\in S_k} \sum_{\substack{\alpha_1\cdots\alpha_{k+1}=\alpha\\ \alpha_1,\alpha_{k+1}\neq\emptyset}} \alpha_1\gamma_{\sigma(1)}'\alpha_2\cdots\gamma_{\sigma(2)}'\cdots \alpha_k \gamma_{\sigma(k)}'\alpha_{k+1}, 
 \end{eqnarray*}
 where in the second equality we used the inductive hypothesis and
 $$\gamma_{\sigma(j)}'= \left\{ \begin{array}{cl}
      \gamma_{\sigma(j)}& \mbox{ if }\sigma(j)\neq i  \\
      \gamma_{i,1}\gamma_{k+1}\gamma_{i,2}& \mbox{ if }\sigma(j)= i
 \end{array} \right. .$$
By subtracting the above equation to the previous one, we obtain
\begin{equation}
    \label{eq:aux2}
    \alpha\{\gamma_1,\ldots,\gamma_{k+1}\} =\sum_{\sigma\in S_k} \sum_{\substack{\alpha_1\cdots \alpha_{k+1} = \alpha \\ \alpha_1,\alpha_{k+1}\neq\emptyset }} \sum_{i=1}^{k+1} \sum_{\substack{\alpha_{i,1}\alpha_{i,2} = \alpha_i\\ \alpha_{1,1},\alpha_{k+1,2}\neq\emptyset}}\alpha_1\gamma_{\sigma(1)}\cdots \alpha_{i,1}\gamma_{k+1}\alpha_{i,2}\cdots \gamma_{\sigma(k)} \alpha_{k+1}.
\end{equation}
Now, given any $\sigma\in S_k$ and $i\in\{1,\ldots,k+1\}$, we define $\tau \in S_{k+1}$ by
$$\tau(j) = \left\{\begin{array}{cl}
     \sigma(j) & \mbox{ if $1\leq j<i$}  \\
     k+1 & \mbox{ if $j=i$}\\
     \sigma(j-1) & \mbox{ if $i<j\leq k+1$}
\end{array} \right. .$$
This map defines a bijection between $S_k\times\{1,\ldots,k+1\}$ and $S_{k+1}$. Hence, Equation \eqref{eq:aux2} can be written as
$$ \alpha\{\gamma_1,\ldots,\gamma_{k+1}\} = \sum_{\tau\in S_{k+1}} \sum_{\substack{ \alpha_1\cdots\alpha_{k+2} = \alpha\\\alpha_1,\alpha_{k+2}\neq\emptyset}} \alpha_1 \gamma_{\tau(1)}\alpha_2\gamma_{\tau(2)} \cdots \alpha_{k+1} \gamma_{\tau(k+1)} \alpha_{k+2},$$
as we wanted to show.
\end{proof}

Recall that $L$ is a  locally finite dimensional connected graded vector space, where the $n$-th graded component is given by the linear span of all the words of length $n$.  The next step is to find a formula for the dual map $\overline{\delta}:=\sum\limits_{n\geq 1}{\delta}_n$ of the previous symmetric brace map on $L$. For this purpose, we introduce the scalar product making $X^*$  an orthonormal basis of $L$:
$$\langle \alpha| w \rangle = \left\{\begin{array}{ll} 1 & \mbox{ if } w = \alpha \\ 0& \mbox{ otherwise}  \end{array} \right. , \qquad\forall\,\alpha , w\in  X^*.$$ The scalar product allows us to identify $L$ with its graded dual $L^\ast$.

Now by duality, the desired coproduct must satisfy that
\begin{equation}
\label{eq:dualeq}
\langle \alpha\{\gamma_1,\ldots,\gamma_n\} | w \rangle = \langle \alpha \otimes \gamma_1\cdots \gamma_n | {\delta}_n(w)\rangle,
\end{equation}
for any $\alpha,\gamma_1,\ldots,\gamma_n\in L$ and $w\in L^*$. Observe that $\gamma_1\cdots \gamma_n$ is a monomial in $\kappa[L]$. Recall that by \eqref{eq:extension}, we know how the scalar product can be extended for products of elements in $L^*$ acting on $\kappa[L]$. 

Notice also that, for degree reasons, ${\delta}_n(w)=0$ if $\operatorname{deg}(w)>n$, so that when $\overline\delta$ acts on a given element $y$ of $L^\ast$, the sum $\sum\limits_{n\geq1}{\delta}_n(y)$ restricts to a finite sum and is therefore well-defined.

\begin{prop}
\label{prop:dualcop}
The map $\overline{\delta}: L^* \to L^* \otimes \kappa[L^*]$ is given by
\begin{equation}
\label{eq:newdualcop}
\overline{\delta}(w) = \sum_{m=1}^\infty \sum_{\substack{w_1\cdots w_{2m+1}=w\\w_1,w_{2m+1},w_{2i}\neq\emptyset,\,\forall\,i}} w_1w_3\cdots w_{2m+1}\otimes w_2\cdot w_4\cdot \cdots \cdot w_{2m},\qquad \forall\;w\in X^*,
\end{equation}
where the product $\cdot$ in the second component is the product in the algebra $\kappa[L^*]$.
\end{prop}
\begin{rem}
Observe that in the right-hand side of \eqref{eq:newdualcop}, the first component of each summand is given by concatenating the subwords $w_1w_3\cdots w_{2m+1}$. On the other hand, the second component of each summand consist of a monomial given by the product in the polynomial algebra $\kappa[L^*]$ of the elements $w_2,w_4\ldots,w_{2m}$. 
\end{rem}

\begin{proof}
 We will check that $\overline{\delta}$ defined in \eqref{eq:newdualcop} is the dual coproduct of the symmetric brace operation, i.e. \eqref{eq:dualeq} holds for any $\alpha,\gamma_1,\ldots,\gamma_n\in L$, $w\in L^*$, $n\geq 1$.
 Indeed, observe that
\begin{eqnarray*}
\langle \alpha\otimes \gamma_1\cdot\cdots\cdot \gamma_n|\overline{\delta}(w) \rangle &=& \sum_{m=1}^\infty \sum_{\substack{w_1\cdots w_{2m+1}=w\\w_1,w_{2m+1},w_{2i}\neq\emptyset,\,\forall\,i}} \langle \alpha| w_1w_3\cdots w_{2m+1} \rangle \langle \gamma_1\cdot\cdots\cdot \gamma_n| w_2\cdot \cdots \cdot w_{2m} \rangle 
\\ &=&  \sum_{\substack{w_1\cdots w_{2n+1}=w\\w_1,w_{2n+1},w_{2i}\neq\emptyset,\,\forall\,i}}  \sum_{\sigma\in S_n }\langle \alpha| w_1w_3\cdots w_{2n+1} \rangle \langle  \gamma_{\sigma(1)} | w_2 \rangle \cdots  \langle  \gamma_{\sigma(n)} | w_{2n} \rangle,
\end{eqnarray*}
where we used \eqref{eq:extension} in the second equality. Fix $\sigma\in S_n$. Observe that the non-vanishing terms in this double sum are such that
\begin{equation}
    \label{eq:cond1}
    \langle \alpha| w_1w_3\cdots w_{2n+1} \rangle \langle  \gamma_{\sigma(1)} | w_2 \rangle \cdots  \langle  \gamma_{\sigma(n)} | w_{2n} \rangle = 1
\end{equation}
which happens for $w = w_1\cdots w_{2n+1}$, $w_1,w_{2n+1},w_{2i}\neq \emptyset$, $\forall$ $i$, when
$\alpha = w_1w_3\cdots w_{2n+1}$ and $w_{2i} = \gamma_{\sigma(i)}$, $\forall$ $i$.

  On the other hand, recalling Proposition \ref{prop:braceproduct} and using the same notation as in the proposition,  we have for the non-vanishing terms arising in the expansion of $\langle \alpha\{\gamma_1,\ldots,\gamma_n\} | w\rangle$:
  
\begin{equation}
\label{eq:cond2}
 \langle \alpha_1\gamma_{\sigma(1)}\cdots \alpha_n\gamma_{\sigma(n)}\alpha_{n+1} |w \rangle = 1 
\end{equation}
which happens when there exist a permutation $\sigma$ and $2n+1$ subwords of $w$, namely $w_1,\ldots,w_{2n+1}$ such that $w = w_1\cdots w_{2n+1},$ 
$ w_{2i-1} = \alpha_i,\, w_{2i} = \gamma_{\sigma(i)}$ for every $1\leq i\leq n$, with $\alpha_1,\alpha_{n+1}\neq \emptyset$ and $\alpha = \alpha_1\cdots \alpha_{n+1}$.

Finally, we notice that terms \eqref{eq:cond1} and \eqref{eq:cond2} are canonically in bijection, and thus obtain the desired conclusion.
\end{proof}

\section{Computing Iterated Pre-Lie Products}
\label{sec:iteratedprelie}

In this section, we will prove some lemmas that will allow us to compute iterated pre-Lie products via the dual coproduct and its iterates; this with the aim of computing later expressions such as the pre-Lie exponential 
\begin{equation}
\label{eq:pre-LieExp}
\exp^\lhd(\alpha) = \alpha + \frac{1}{2!}\alpha\lhd \alpha + \frac{1}{3!} (\alpha\lhd \alpha) \lhd \alpha+\cdots 
\end{equation}
of an element of a pre-Lie algebra $(L,\lhd)$. 

\begin{nota}
Let  $(L,\lhd)$ be a connected graded and locally finite dimensional pre-Lie algebra. 
\begin{enumerate}
    \item We write again $\overline\delta=\sum\limits_{n\geq 1}{\delta}_n$ for the dual map of the symmetric brace product given by Theorem \ref{thm:OudomGuin}.
    \item Recall that the Hopf algebra $(\kappa[L],\ast,\Delta)$ is an enveloping algebra for $L$. We write $(\kappa[L^*],\cdot,\delta)$ for the graded dual Hopf algebra. The product $\cdot$ is simply the product of polynomials.
    \item We set $\delta_{\mathrm{irr}}=\delta_{\mathrm{irr}}^{[2]}:={\delta}_1$. In general, for $n\geq2$, $\delta^{[n]}_{\mathrm{irr}}$ stands for the restriction of the iterated coproduct ${\delta}^{[n]}$ from $\kappa[L^*]$ to $\kappa[L^*]^{\otimes n}$ to a map from $L^*$ to $(L^*)^{\otimes n}$ (by restriction to $L^\ast\subset \kappa[L^*]$ and composition with the canonical projection from $\kappa[L^*]^{\otimes n}$ to $(L^*)^{\otimes n}$). That ${\delta}_1$ identifies with the restriction of $\delta$ to a map from $L^*$ to $(L^*)^{\otimes 2}$ follows from \cite[Lemma 6.2.1]{CP21}.
    \item Let $ \alpha,\beta\in L$ and $w\in L^*$. By considering $\alpha$ as an element in $({L^*})^{*}$, we define $\hat{\alpha}$ to be the infinitesimal character on $\kappa[L^*]$ such that $\hat{\alpha}$ coincides with $\alpha$ on $L^*$. By infinitesimal character, we mean that for $w_1,\ldots,w_n\in L^*$, $\langle \hat\alpha | w_1\cdots w_n \rangle = 0$ whenever $n\neq 1$.
\end{enumerate}
\end{nota}

\begin{lemma}
\label{lem:iterated1}
Let $(L,\lhd)$ be a connected graded and locally finite dimensional pre-Lie algebra. Then
$$ \langle \alpha \lhd \beta | w\rangle = \langle \alpha\otimes\beta | \delta_{\mathrm{irr}}(w)\rangle = \langle \hat{\alpha}\otimes \hat{\beta} | \overline{\delta}(w) \rangle,$$
for any  $\alpha,\beta \in L$ and $w\in L^*$, where $\overline{\delta}(x) := \delta(x) - 1\otimes x - x\otimes 1$ is the reduced coproduct on $\kappa[L^*]$.
\end{lemma}
\begin{proof}
Let $\alpha,\beta\in L$ and $w\in L^*$. Using duality, we have that
\begin{eqnarray*}
\langle \alpha\otimes\beta | \delta_{\mathrm{irr}}(w)\rangle &=& \langle \alpha\otimes\beta | {\delta}_1(w)\rangle
\\ &=& \langle \alpha\{\beta\} | w\rangle
\\ &=& \langle \alpha \lhd \beta | w \rangle ;
\end{eqnarray*} 
this proves the first equality. For the second equality, by \cite[Lemma 6.2.1]{CP21}, one has for the restriction of $\overline{\delta}$ to $L^*$: $$\overline{\delta} = \delta_{\mathrm{irr}} + \sum\limits_{n\geq 2}\delta_n.$$ 
Since ${\delta}_n(w) \in L^* \otimes \kappa[L^*]_n$, we get that
\begin{eqnarray*}
\langle  \hat{\alpha}\otimes \hat{\beta} | \overline{\delta}(w) \rangle &=& \langle  \hat{\alpha}\otimes \hat{\beta} | \delta_{\mathrm{irr}}(w) +\sum_{n\geq 2}\delta_n(w)\rangle
\\&=&  \langle  \hat{\alpha}\otimes \hat{\beta} | \delta_{\mathrm{irr}}(w)\rangle + 0
\\ &=& \langle  {\alpha}\otimes {\beta} | \delta_{\mathrm{irr}}(w)\rangle,
\end{eqnarray*}
where in the last equality, we used that $\hat\alpha$ and $\hat\beta$ coincide with $\alpha$ and $\beta$ on $L^*$, respectively. 
\end{proof}
One can obtain the following generalization of the previous result. This will allow us to compute the iterated pre-Lie products via the iteration of the restricted coproduct.

\begin{lemma}
\label{lem:iterated2}
Let $(L,\lhd)$ be a connected graded and locally finite dimensional pre-Lie algebra. Then
\begin{equation}
    \label{eq:iterated}
    \langle (\cdots (\alpha_1\lhd \alpha_2)\lhd \cdots )\lhd \alpha_n |w \rangle = \langle \alpha_1\otimes \cdots\otimes \alpha_n | \delta^{[n]}_{\mathrm{irr}}(w)\rangle = \langle \hat{\alpha}_1\otimes \cdots \otimes \hat{\alpha}_n | \bar{\delta}^{[n]}(w) \rangle,
\end{equation}

for any $n\geq 2$, $\alpha_1,\ldots,\alpha_n\in L$ and $w\in L^*$.
\end{lemma}
\begin{proof}  We prove the lemma by induction on $n$. The case $n=2$ is given in Lemma \ref{lem:iterated1}. Now, assume that the lemma holds for $n-1$ and define $\gamma: =(\cdots (\alpha_1\lhd \alpha_2)\lhd \cdots )\lhd \alpha_{n-1}\in L.$ Also, introduce the Sweedler-type notation $\delta_{\mathrm{irr}}(w) = w^{\prime(1)}\otimes w^{\prime(2)} $. Hence
\begin{eqnarray*}
\langle \gamma \lhd \alpha_n|w\rangle &=& \langle \gamma\otimes \alpha_n |   w^{\prime(1)}\otimes w^{\prime(2)} \rangle
\\ &=& \langle \gamma | w^{'(1)}\rangle \langle \alpha_n| w^{'(2)}\rangle
\\ &=& \langle \alpha_1\otimes \cdots \otimes \alpha_{n-1} | \delta_{\mathrm{irr}}^{[n-1]}(w^{'(1)}) \rangle \langle \alpha_{n}| w^{\prime(2)}\rangle 
\\ &=& \langle \alpha_1\otimes \cdots\otimes \alpha_n | \delta_{\mathrm{irr}}^{[n-1]}(w^{'(1)}) \otimes w^{'(2)}\rangle
\\ &=& \langle \alpha_1 \otimes \cdots \otimes \alpha_n | \delta_{\mathrm{irr}}^{[n]}(w)\rangle,
\end{eqnarray*}
where  we used the inductive hypothesis in the third equality and that $ \delta_{\mathrm{irr}}^{[n]} = (\delta_\mathrm{irr}^{[n-1]}\otimes \Id) \circ \delta_\mathrm{irr}$ in the last equality. This identity, in turn, follows from
\begin{enumerate}
    \item the coassociativity of $\delta$ that implies $ \delta^{[n]} = (\delta^{[n-1]}\otimes \Id) \circ \delta$,
    \item the definition of $\delta_{\mathrm{irr}}^{[n]}$ by restriction of $\delta^{[n]}$ to a map from $L^*$ to $(L^*)^{\otimes n}$.
\end{enumerate}

The proof of the second equality in \eqref{eq:iterated} follows e.g. by induction from Lemma \ref{lem:iterated1} and the coassociativity of the reduced coproduct $\overline{\delta}$ (see e.g. \cite[Def. 2.3.5]{CP21}).
\end{proof}

Finally, the following lemma establishes the analogous result to compute  iterated symmetric brace products via the reduced coproduct.

\begin{lemma}
\label{lem:iteratedbrace}
Let $(L,\lhd)$ be a connected graded and locally finite dimensional pre-Lie algebra. Then 
\begin{eqnarray}
   \nonumber
    \Big\langle (\cdots (\alpha_{1,1}\{\alpha_{1,2}&\hspace{-0.8cm},&\hspace{-0.8cm}\ldots\alpha_{m_2,2}\})\cdots )\{\alpha_{1,n},\ldots,\alpha_{m_n,n}\} \;|\: w\Big\rangle \\ &=& \Big\langle \alpha_{1,1}\otimes \alpha_{1,2}\cdots \alpha_{m_2,2}\otimes \cdots \otimes \alpha_{1,n}\cdots \alpha_{m_n,n}\;|\; \overline\delta^{[n]}(w)\Big\rangle  \label{eq:iteratedbrace}
\end{eqnarray}

for any $n\geq 2$, $\alpha_{i,j}\in L$ for $1\leq i\leq n$, $1\leq j\leq n_i$ $(n_1=1)$, and $w\in L^*$.
\end{lemma}

Indeed, $(\dots(\alpha_{1,1}\{\alpha_{1,2},\ldots,\alpha_{m_2,2}\})\cdots )\{\alpha_{1,n},\ldots,\alpha_{m_n,n}\}$ is the projection on $L$ orthogonally to the other graded components of $\kappa[L]$ of the element 
$\alpha_{1,1}\ast(\alpha_{1,2}\ldots\alpha_{m_2,2})\ast\dots\ast (\alpha_{1,n},\ldots,\alpha_{m_n,n})$ in $\kappa[L]$, where $\ast$ is the associative product defined in \eqref{eq:product}. The lemma follows as the $\ast$ product in $\kappa[L]$ is dual to the coproduct $\delta$ in $\kappa[L^\ast]$.

\section{Forest Formulas for Iterated Coproducts}
\label{sec:forestformulas}

We now introduce the main tool for evaluating the pre-Lie exponential: the forest-type formulas for iterated reduced coproducts introduced in \cite{mp}. In this section, we will describe the required notation and results in order to state the indicated formulas.  We follow closely the arguments and notation in that article, in particular we follow almost verbatim the explanation of the notation by means of graphical arguments, to make the present article self-contained. We also introduce symmetry factors, give a different proof of the main identity and extend the forest-type formula for iterated reduced coproducts to iterated coproducts and iterated irreducible coproducts.

Again, let $(L,\lhd)$ be a connected graded and locally finite dimensional pre-Lie algebra. We consider the Hopf algebra $(\kappa[L^\ast],\cdot,\delta)$ and assume that  $L^\ast$ has a basis $\mathcal B=\{ b_i\}_{i\in \mathbb{N}^*}$. Notice that our arguments work also with a finite basis, we make the assumption that the basis is countable only for notational convenience, as in the examples we will consider later this will be the case.

We will also consider multisets over $\mathbb{N}^*$ and, for any such multisets $I,J$, write $I\cup J$ for their union. For example, if $I = \{1,2,2\}$ and $J=\{2,3,3\}$, then $I\cup J=\{1,2,2,2,3,3\}$. Hence, one can take monomials $b_I=\prod_{i\in I}b_i$, so that $b_I\cdot b_J=b_{I\cup J}$. Observe that $b_\emptyset$ can be taken as the unit of $\kappa[L^*]$.

In the Hopf algebra $\kappa[L^*]$, the reduced coproduct of an element $b_i\in \mathcal{B}$ can be expanded as
\begin{equation}
    \label{eq:generalforest0}
    \overline\delta (b_i)=\sum_{i_0, I\not=\emptyset} \lambda^{i ; i_0}_{I} b_{i_0}\otimes b_{I},
\end{equation}
for coefficients $\lambda^{i;i_0}_I\in \kappa$, with $i_0\in \mathbb{N}^*$ and $I\subseteq \mathbb{N}^*$ a non-empty multiset. These coefficients completely determine the coproduct and its action on products, as well as the action of the iterated coproducts. 
\par For our purposes, it will be beneficial to index the above summation by non-planar decorated corollas in the following way:

\begin{equation}
\label{eq:generalforest}
\overline\delta (b_i)=\sum \lambda\Big(\corollaa \Big) b_{i_0}\otimes b_{i_1}\dots b_{i_k},
\end{equation}
where the decoration is given by a pair of integers $(i;i_0)$ associated to the root, and leaves are decorated by positive integers $i_1,\ldots,i_k$. Here, the non-planarity of the rooted trees means that the ordering of the branches does not matter, reflecting the commutativity of the product in the Hopf algebra.

\begin{ex}
Consider a single term in the sum \eqref{eq:generalforest} for the reduced coproduct such that the indexing multiset is of the form $I = \{i_1,i_2\}$, namely
\begin{equation}
    \label{eq:I2term}
    \lambda\Big (\corollab \Big ) b_{i_0}\otimes b_{i_1} b_{i_2}.
\end{equation}
By definition of reduced coproduct, we have
\[
\overline\delta(b_{i_1} b_{i_2}) = (1\otimes b_{i_1}+b_{i_1} \otimes 1+\overline\delta(b_{i_1}))(1\otimes b_{i_2} +b_{i_2}\otimes 1+\overline\delta(b_{i_2}))-1\otimes b_{i_1}b_{i_2}-b_{i_1}b_{i_2}\otimes 1.
\]
Observe that the contribution of \eqref{eq:I2term} to $\overline\delta^{[3]}(b_i)=(\Id\otimes \overline\delta)(\overline\delta(b_i))$ will split in four terms, whose complexity is encoded by the appearance of products of coefficients $\lambda^{j;j_0}_{J}$:
\begin{enumerate}
    \item There is a first term with no more complexity than in $\overline\delta(b_i)$:
\[\lambda
\Big (\corollab \Big ) b_{i_0}\otimes (b_{i_1}\otimes b_{i_2} +b_{i_2}\otimes  b_{i_1}).
\]
\item There is a second term where only the reduced coproduct of $b_{i_1}$ occurs:  
\[\lambda\Big (\corollab \Big ) \left(\sum \lambda^{i_1;i_{1,0}}_{i_{1,1},\dots ,i_{1,k}} b_{i_0}\otimes (b_{i_{1,0}}\otimes b_{i_2}b_{\{i_{1,1},\dots, i_{1,k}\}} +b_{i_2}b_{i_{1,0}}\otimes  b_{\{i_{1,1},\dots, i_{1,k}\}})\right).
\] This term arises from the product $(1\otimes b_{i_2}+ b_{i_2}\otimes 1)\overline{\delta}(b_{i_1})$ and is naturally indexed by the trees
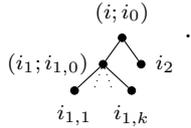
\[\corollac .\]
\item In the same way there is a contribution, corresponding to $(1\otimes b_{i_1}+b_{i_1} \otimes 1)\overline\delta(b_{i_2})$, indexed by the trees
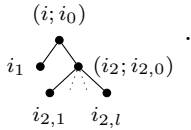
\[\corollad .\] 
\item Finally the terms in relation with $\overline\delta(b_{i_1})\overline\delta(b_{i_2})$ that will be indexed by trees 
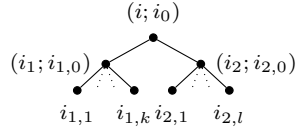
\[\corollae .\]
\end{enumerate}
\end{ex}

When iterating the reduced coproduct, such groups of terms naturally appear, labelled by trees that encode the presence of the coefficients $\lambda^{j;j_0}_{j_1,\ldots,j_k}$.

In order to state the forest formula for the iterated reduced coproduct, we will need the subsequent notions associated to decorated rooted trees.

\begin{defi}\label{def:tree}
Let $T$ be a finite rooted tree whose internal vertices are decorated by pairs $p=(p_1;p_2)$ of positive integers, and leaves are decorated by positive integers. In the case of $T$ being a  single-vertex tree, the vertex is considered as a leaf.
\begin{enumerate}
    \item Let $V(T)$ be the set of vertices of $T$. In addition, let $\operatorname{Int}(T)$ and $\operatorname{Leaf}(T)$ be the sets of internal vertices of $T$ and leaves of $T$, respectively. For any $x\in \operatorname{Int}(T)$, we denote $d(x)=(d_1(x);d_2(x))$ its decoration, and, if $x\in \operatorname{Leaf}(T)$, we denote for convenience its decoration $d(x)=d_1(x)=d_2(x)$. For any internal vertex $x$, we also note, $\operatorname{succ}(x)$ the set of its immediate successors.
    \item If the root of $T$ is decorated by $i$ or $(i;i_0)$, we say that the tree is \textit{associated to $b_i$}. We will denote it by $T\in \mathcal{T}_{i}$.
    \item For a given pair $p$, we denote by $B^+_p(T_1,\ldots, T_s)$ the tree obtained by adding a common root decorated by $p$ to the trees $T_1,\ldots, T_s$. {If $T = B^+_p(T_1,\ldots, T_s)$, we will denote by $B^-(T)$ the multiset of trees $T_1\cdots T_s$.}
    \item We define the length of $T$, denoted by $|T|$, as the number of elements in $T$ viewed as a poset. 
    \item Let $F$ be a multiset of decorated trees (that is, repetitions are allowed). We write $F$ as 
    \begin{equation}
        \label{eq:multiforest}
        F=\{T_{1,i_1}^{k_{1,1}},\dots,T_{s_1,i_1}^{k_{s_1,1}}\}\cup \dots \cup \{T_{1,i_p}^{k_{1,p}},\dots,T_{s_p,i_p}^{k_{s_p,p}}\}
    \end{equation}
    where: \begin{itemize}
        \item the tree $T_{j,i_q}$ is associated to $i_q$;
        \item the trees $T_{j,i_q}$ are all distinct;
        \item the notation $T_{j,i_p}^{k_{j,p}}$ means that the tree $T_{j,i_p}$ appears with multiplicity $k_{j,p}$ in the multiset $F$.
    \end{itemize}  
    
    The {\it symmetry coefficient} of $F$, $\operatorname{sym}(F)$ is then, by definition $$\operatorname{sym}(F):=\prod_{j=1}^p {{k_{1,j}+\cdots+k_{s_j,j}}\choose{k_{1,j},\dots,k_{s_j,j}}}.$$
    \item We define the coefficient $\lambda(T)$ as follows: if $\bullet_i$ stands for the single-vertex tree with decoration $i$, then $\lambda(\bullet_i):=1$. More generally, if $T=B^+_{(i;i_0)} (T_1,\ldots, T_s)$, then 
\begin{equation}
    \label{def:lambdaT}
    \lambda(T):=\lambda^{i;i_0}_{i_1,\dots,i_s}\cdot \operatorname{sym}(F)\cdot \lambda(T_1)\cdots\lambda(T_s)
\end{equation}
when $T_1,\dots, T_s$ are trees respectively associated to $b_{i_1},\ldots,b_{i_s}$, and where $F$ is the multiset $\{T_1,\ldots,T_s\}$. In other words, if for $x\in V(T)$, $T^x$ stands for the subtree of $T$ consisting of $x$ and all its descendants, then
{
\begin{equation}
    \lambda(T) = \prod_{x\in \operatorname{Int}(T)} \lambda^{d(x)}_{d_1(\operatorname{succ}(x))} \operatorname{sym}(B^-(T^x)),
\end{equation}}
where $d_1(\operatorname{succ}(x)) = \{d_1(y)\,:\, y\in \operatorname{succ}(x) \}$.

\end{enumerate}
\end{defi}

In \cite{mp}, it was observed that, since the decorations of trees associated to $b_i$ run over all the indexes of basis elements appearing in the various iterated reduced coproducts of $b_i$, a way the describe all the tensors of length $k$ that can be obtained in the $k$-fold iterated reduced coproduct is by using the notion of $k$-linearization of a poset.

\begin{defi}
\label{def:Cf}
Let $P$ be a finite poset of cardinality $n$. A \textit{linearization} of $P$ is a bijective, strictly order preserving map $f:P\to[n]$. A \textit{$k$-linearization} of $P$ is a surjective, strictly order preserving map $f:P\to[k]$. A \textit{weak $k$-linearization} of $P$ is a strictly order preserving map $f:P\to[k]$. We denote by 
$\mathrm{lin} ( P)$, resp. $k-\mathrm{lin} ( P)$, resp. $\mathrm{w}-k-\mathrm{lin} ( P)$, the set of linearizations, resp. $k$-linearizations of $P$, resp. weak $k$-linearizations of $P$. 
\end{defi}

A $k$-linearization is a surjective weak $k$-linearization; a linearization a bijective $k$-linearization.

Notice, for later use, that if $P$ is a disjoint union of posets $P_1\sqcup\cdots\sqcup P_q$, with no relations between the elements of two distinct components, then weak $k$-linearizations of $P$ are in bijection with $p${-tuples} of weak $k$-linearizations of the components $P_1,\dots,P_q$. This observation does not hold for $k$-linearizations.

\begin{defi}
Let $T$ be a tree with a given decoration $d=(d_1;d_2)$. If $f$ is a linearization of $T$, a $k$-linearization of $T$, or a weak $k$-linearization of $T$, we write
$$C(f):=\left(\prod_{x_1\in f^{-1}(\{1\})} b_{x_1}\right) \otimes \cdots \otimes \left(\prod_{x_k\in f^{-1}(\{k\})} b_{x_k}\right),$$
where $b_{x_i}$ stands for $b_{d_2(x_i)}$.

In the formula, when $f^{-1}(\{i\})=\emptyset$, , we set $\left(\prod_{x_i\in f^{-1}(\{i\})} b_{x_i}\right):=1$. This can happen only with weak $k$-linearizations.
\end{defi}

\begin{theo}[Forest formulas for the iterated coproduct and the iterated reduced coproduct]
\label{thm:forestform0}
For any $b_i \in \mathcal{B}$, we have for the action of the $k$-fold iterated reduced coproduct, respectively the $k$-fold iterated coproduct, on $\kappa[L^*]$:
\begin{equation}
\label{eq:forestformula0}
\overline\delta^{[k]}(b_i)=\sum\limits_{T\in \mathcal{T}_i}\sum\limits_{f\in k-\mathrm{lin} (T)}\lambda(T)C(f).
\end{equation}
\begin{equation}
\label{eq:forestformulafull}
\delta^{[k]}(b_i)=\sum\limits_{T\in \mathcal{T}_i}\sum\limits_{f\in \mathrm{w}-k-\mathrm{lin} (T)}\lambda(T)C(f).
\end{equation}
\end{theo}

The first part of the theorem was obtained in {\cite[Lemma 12]{mp}}, but the calculation of $\lambda(T)$ did not take into account symmetry factors. We propose here a different proof. 

\begin{proof}
Let us start with some general preliminaries on iterated coproducts.
Let $(C,\Delta,\varepsilon,\eta)$ be a coaugmented coalgebra over a ground field $\kappa$. The coproduct $\Delta$ is assumed to be coassociative and counital, with counit $\varepsilon: C\to \kappa$. The coaugmentation $\eta$ maps $\kappa$ to $C$. One sets $\nu:=\eta\circ\varepsilon$. The iterated coproduct and the iterated reduced coproduct are then related by the identity:
$$\overline{\Delta}^{[k]}=(\Id-\nu)^{\otimes k}\circ \Delta^{[k]}.$$
To conversely express $\Delta^{[k]}$ in terms of iterated reduced coproducts, notice that
$$\Delta^{[k]}=(\nu+(\Id-\nu))^{\otimes k}\circ \Delta^{[k]}$$
$$=\left(\sum\limits_{l\leq k}\sum\limits_{1\leq i_1<\dots<i_l\leq k}(\Id-\nu)\otimes\dots\otimes \nu\otimes\dots\otimes\nu\otimes\dots\otimes (\Id-\nu)\right)\circ \Delta^{[k]},$$
where the term in the last summation formula contains $l$ copies of $\nu$, in positions $i_1,\dots,i_l$ (with possibly $i_1=1$ and/or $i_l=k$, that is, $\nu$ is allowed to be in first or last position in the tensor product).
Let now $f:[k]\hookrightarrow[n]$ be an increasing injection. We denote $\hat{f}$ the map from  $C^{\otimes k}$ to $C^{\otimes n}$ defined by 
$$\hat{f}(c_1\otimes\dots\otimes c_k):=d_1\otimes\dots\otimes d_n$$
with $d_{f(i)}:=c_i$ for $1\leq i\leq k$, and $d_j:=1$ if $j$ is not in the image of $f$.
Since $\Delta$ is counital, the above expression of $\Delta^{[k]}$ rewrites 
$$\sum\limits_{l\leq k}\sum\limits_{f:[k-l]\hookrightarrow[k]}\hat{f}\circ \overline{\Delta}^{[k-l]}.$$

Now, a weak $k$-linearization can be uniquely obtained as the composition of a $l$-linearization, with $l\leq k$, and an injection $f$ from $[l]$ into $[k]$. This is a particular case of the usual decomposition of a map between two finite sets as a surjection followed by an injection. The equivalence of (\ref{eq:forestformula0}) and (\ref{eq:forestformulafull}) follows.

Let us prove these formulas jointly by induction on $k$. 
We begin with $$\overline{\delta}(b_i)=\sum\limits_{i_0,i_1,\dots,i_n}\lambda^{i;i_0}_{i_1,\ldots,i_n}b_{i_0}\otimes b_{i_1}\dots b_{i_n}.$$
Formula (\ref{eq:forestformula0}) is then obvious for $k=2$ since it is a rewriting in graphical terms of this expansion as in \eqref{eq:generalforest}. Let us assume $k>2$. Then we get
$$\overline{\delta}^{[k]}(b_i)=\sum\limits_{i_0,i_1,\dots,i_n}\lambda^{i;i_0}_{i_1,\dots,i_n}b_{i_0}\otimes \overline{\delta}^{[k-1]}(b_{i_1}\cdots b_{i_n}).$$
We are left with the problem of computing $$\overline{\delta}^{[k-1]}(b_{i_1}\cdots b_{i_n})=(\Id-\nu)^{\otimes k-1}\circ {\delta}^{[k-1]}(b_{i_1}\cdots b_{i_n}).$$
However, since the coproduct is a morphism of algebras, we have that
$${\delta}^{[k-1]}(b_{i_1}\cdots b_{i_n})={\delta}^{[k-1]}(b_{i_1})\cdots {\delta}^{[k-1]}(b_{i_n}).$$
By the induction hypothesis,
$$\delta^{[k-1]}(b_{i_j})=\sum\limits_{T\in \mathcal{T}_{i_j}}\sum\limits_{f\in \mathrm{w}-{(k-1)}-\mathrm{lin} (T)}\lambda(T)C(f)$$
and
\begin{eqnarray}
{\delta}^{[k-1]}(b_{i_1})\cdots {\delta}^{[k-1]}(b_{i_n})
&=&\prod\limits_{j=1}^n\left(\sum\limits_{T\in \mathcal{T}_{i_j}}\sum\limits_{f\in \mathrm{w}-{(k-1)}-\mathrm{lin}\nonumber (T)}\lambda(T)C(f)\right)\\
&=&\sum\limits_{\substack{T_j\in \mathcal{T}_{i_j}\\ 1\leq j\leq n}}\lambda'(T)\left(\sum\limits_{f\in {(k-1)}-\mathrm{lin} (T)}C(f)\right),\label{eq:auxforest1}
\end{eqnarray}
where in the last equality, $T$ is the multiset of decorated trees $T_1\cdots T_n$. That is, as a poset, $T$ is the disjoint union of the $T_i$. Also, we have set $$\lambda'(T):=\lambda(T_1)\cdots \lambda(T_n).$$
Moreover, in \eqref{eq:auxforest1} we used the fact that weak $(k-1)$-linearizations of $T$ are in bijection with families of $(k-1)$-linearizations of the posets $T_1,\ldots,T_n$ in order to obtain 
\begin{eqnarray}
\overline{\delta}^{[k-1]}(b_{i_1}\cdots b_{i_n})&=&(\Id-\nu)^{\otimes k-1}\circ \delta^{[k-1]}(b_{i_1}\cdots b_{i_n})\nonumber
\\&=&\sum\limits_{\substack{T_j\in \mathcal{T}_{i_j}\\ 1\leq j\leq n}}\lambda'(T)\left(\sum\limits_{f\in {(k-1)}-\mathrm{lin} (T)}C(f)\right). \label{eq:auxforest}
\end{eqnarray}

Let us say $T\in \mathcal{T}_{i_1,\dots,i_n}$ when $T$ is the forest $T_1\cdots T_n$ and there exists a permutation $\sigma$ of $[n]$ such that each tree $T_j $ is associated to a $b_{i_{\sigma(j)}}$. Hence we can rewrite \eqref{eq:auxforest} as
\begin{equation}
    \label{eq:auxforest2}
    \overline{\delta}^{[k-1]}(b_{i_1}\cdots b_{i_n})=\sum\limits_{T\in \mathcal{T}_{i_1,\ldots,i_n}}\lambda'(T)\gamma(T)\left(\sum\limits_{f\in {(k-1)}-\mathrm{lin} (T)}C(f)\right),
\end{equation}
for a certain coefficient $\gamma(T)$ to be determined.

Let us rewrite now $b_{i_1}\cdots b_{i_n}$ as a monomial $b_{j_1}^{p_1}\cdots b_{j_m}^{p_m}$, with the $b_{j_h}$ pairwise distinct.
Viewing now $T$ as a multiset of trees, we use the same notation as in Definition \ref{def:tree}.5 and express $T$ as $$T=\{T_{1,j_1}^{k_{1,1}},\dots,T_{s_1,j_1}^{k_{s_1,1}}\}\cup \cdots \cup \{T_{1,j_m}^{k_{1,m}},\dots,T_{s_m,j_m}^{k_{s_m,m}}\}.$$
We then have $k_{1,q}+\cdots +k_{s_q,q}=p_q$ and $p_1+\cdots +p_m=n$. The multiplicity of $T$ in the right-hand side of the expansion of $ \overline{\delta}^{[k-1]}(b_{i_1}\cdots b_{i_n})$ in \eqref{eq:auxforest2} is then obtained as the product over $1\leq j\leq m$ of the number of ordered partitions of a set $X$ of cardinal $p_j$ into a disjoint union $X_1\sqcup\cdots\sqcup X_{s_j}$ with $|X_l|=k_{l,j}$ for $1\leq l\leq s_j$. This is precisely
the symmetry coefficient of $T$, $\operatorname{sym}(T)$, that is $\gamma(T) = \operatorname{sym}(T)$ and hence
$$\overline{\delta}^{[k-1]}(b_{i_1}\cdots b_{i_n})=\sum\limits_{T\in \mathcal{T}_{i_1,\ldots,i_n}}\lambda(T)\operatorname{sym}(T)\left(\sum\limits_{f\in {(k-1)}-\mathrm{lin} (T)}C(f)\right),$$
Therefore (\ref{eq:forestformula0}) and the equivalent formula (\ref{eq:forestformulafull}) follow as we wanted to prove.
\end{proof}

The above theorem leads us to the following forest-type formula for the iterated restricted coproduct $\delta_{\mathrm{irr}}^{[n]}$ on $L^*$.
\begin{theo}[Forest formula for the irreducible coproduct]
For any $b_i\in \mathcal{B}$, we have for the action of $\delta_{\mathrm{irr}}^{[k]}$:
\begin{equation}
\label{eq:forestformula}
\delta_{\mathrm{irr}}^{[k]}(b_i)=\sum\limits_{T\in \mathcal{T}_i}\sum\limits_{f\in \mathrm{lin} (T)}\lambda(T)C(f).
\end{equation}
\end{theo}
\begin{proof}
Recall that $\delta_\mathrm{irr}^{[k]}$ is the restriction of the iterated coproduct $\delta^{[k]}:\kappa[L^*]\to \kappa[L^*]^{\otimes k}$ to a map $L\to  (L^*)^{\otimes k}$. That is, $\delta^{[k]}(b_i)$ can be split into two terms, namely ${\delta}_\mathrm{irr}^{[k]}(b_i)$ and a linear combination of terms $b_{I_1}\otimes\cdots \otimes b_{I_k}$ such that at least one of the $b_{I_j}$ is equal to the unit 1 or is such that $|I_j|\geq 2$. This second term is projected to 0 when we restrict $\delta^{[k]}$ in order to obtain $\delta_\mathrm{irr}^{[k]}$. 
\par Now, consider a term $C(f)$ in the forest formula \eqref{eq:forestformula0} for $\overline{\delta}$ given in Theorem \ref{thm:forestform0}. Since $f$ is a surjective map, then $l(T)\geq k$. In the case that $f$ is not a bijection, i.e. $l(T)>k$, there exists $j\in [k]$ such that $f^{-1}(j)$ contains at least two elements. In this case we would have that $\prod_{x_j\in f^{-1}(j)} b_{x_j}\not\in L^*$ and thus $C(f)$ is projected to 0 when we consider $\delta_\mathrm{irr}^{[k]}$.
\par On the other hand, take $C(f)$ in \eqref{eq:forestformula0} such that $f$ is a bijection. It easily follows that $C(f) \in (L^*)^{\otimes k}$ and hence this term appears in $\delta_\mathrm{irr}^{[k]}(b_i)$. The conclusion is that $\delta_\mathrm{irr}^{[k]}$ also has a forest-type formula obtained from \eqref{eq:forestformula0} by only considering $k$-linearizations of $T$ that are bijections, i.e. linearizations of $T$. We obtain \eqref{eq:forestformula}.
\end{proof}

\section{Pre-Lie Exponential in the Free Pre-Lie Algebra}
\label{sec:preliexexpfree}

We are interested in computing the action of the pre-Lie exponential (or Agrachev-Gamkrelidze operator) on the generator of the free pre-Lie algebra of rooted trees. The results are classical and relate to Runge-Kutta methods \cite{Brouder,Butcher}, but we use this example to show how the machinery of the present article works in practice on this simple case.
    \par In order to state the definition of the pre-Lie algebra of rooted trees, we consider the countable set of non-planar rooted trees
    $$\mathcal{T} = \left\{ \begin{tikzpicture}[scale=0.5,
level 1/.style={level distance=7mm,sibling distance=7mm}]
\node(0)[solid node,label=above:{}]{};
\end{tikzpicture}
 , \begin{tikzpicture}[baseline={([yshift=-1.4ex]current bounding box.center)},scale=0.5,
level 1/.style={level distance=7mm,sibling distance=7mm}]
\node(0)[solid node,label=above:{}]{}
child{node(1)[solid node]{}};
\end{tikzpicture} ,  \begin{tikzpicture}[baseline={([yshift=-1.4ex]current bounding box.center)},scale=0.5,
level 1/.style={level distance=7mm,sibling distance=7mm}]
\node(0)[solid node,label=above:{}]{}
child{node(1)[solid node]{}
child{[black] node(11)[solid node]{}}
};
\end{tikzpicture} , \begin{tikzpicture}[baseline={([yshift=-1.4ex]current bounding box.center)},scale=0.5,
level 1/.style={level distance=7mm,sibling distance=7mm}]
\node(0)[solid node,label=above:{}]{} 
child{node(1)[solid node]{}
}
child{node(2)[solid node]{}
};
\end{tikzpicture} , \begin{tikzpicture}[baseline={([yshift=-1.4ex]current bounding box.center)},scale=0.5,
level 1/.style={level distance=7mm,sibling distance=7mm}]
\node(0)[solid node,label=above:{}]{}
child{node(1)[solid node]{}
child{[black] node(11)[solid node]{}
child{[black] node(111)[solid node]{}}
}
};
\end{tikzpicture} , \begin{tikzpicture}[baseline={([yshift=-1.4ex]current bounding box.center)},scale=0.5,
level 1/.style={level distance=7mm,sibling distance=7mm}]
\node(0)[solid node,label=above:{}]{}
child{node(1)[solid node]{}
child{[black] node(11)[solid node]{}}
child{[black] node(111)[solid node]{}}
};
\end{tikzpicture} , \begin{tikzpicture}[baseline={([yshift=-1.4ex]current bounding box.center)},scale=0.5,
level 1/.style={level distance=7mm,sibling distance=7mm}]
\node(0)[solid node,label=above:{}]{} 
child{node(1)[solid node]{}
child{[black] node(11)[solid node]{}}
}
child{node(2)[solid node]{}
};
\end{tikzpicture}
,   \begin{tikzpicture}[baseline={([yshift=-1.4ex]current bounding box.center)},scale=0.5,
level 1/.style={level distance=7mm,sibling distance=7mm}]
\node(0)[solid node,label=above:{}]{} 
child{node(1)[solid node]{}
}
child{node(2)[solid node]{}
}
child{[black] node(11)[solid node]{}
};
\end{tikzpicture},\ldots\right\},$$
For later use, we order the elements of $\mathcal{T}$ to put them in bijection with the set of positive integers $\mathbb{N}$, i.e. {$\T=\{t_j\}_{j\in\mathbb{N}^*}$.} We also assume that $t_0=\bullet$ is the single-vertex tree.

For $t,t'\in \mathcal{T}$, we define the element $t\lhd t'$ as the sum of the trees obtained by grafting the root of $t'$ to the vertices of $t$:
\begin{equation}
    \label{eq:pre-LieTrees}
    t\lhd t'= \sum_{v\in V(t)} t \leftarrow_v t',
\end{equation}
where $V(t)$ stands for the vertex set of $t$ and $t \leftarrow_v t'$ stands for the grafting of the root of $t'$ via a new edge to vertex $v$ of $t$. This operation can be linearly extended to the linear span of $\mathcal{T}$, which we denote $L_{\mathrm{free}}$. See e.g. \cite[Section 6.3]{CP21} for details.

    \begin{theo}[\cite{CHL}]
    $(L_{\mathrm{free}},\lhd)$ is a free pre-Lie algebra over the generator $\bullet$ (the single-vertex tree).
    \end{theo}
    
    It is clear that $\Lfr$ is a connected graded and locally finite dimensional pre-Lie algebra, where the grading is given by $\deg(t) = |t|$ with $|t|$ being the number of vertices of $t$. We can therefore apply the formula \eqref{eq:forestformula} together with Lemma \ref{lem:iterated2} to obtain an expression for iterated pre-Lie products.
    \par Let us describe the ingredients for the forest formula. The algebra $\kappa[\Lfr]$, viewed as a Hopf algebra and equipped with the basis of forests of non-planar rooted trees, is called the \textit{Grossman-Larson Hopf algebra}. Here, a forest is a multiset of trees, or equivalently a monomial of trees. The associative product is given by the formula \eqref{eq:product} and can be described in terms of grafting operations.
    \par The dual Hopf algebra $\kappa[\Lfr^*]$ identifies with $\kappa[\Lfr]$ as a vector space. It is a polynomial algebra over non-planar rooted trees. It is often called the \textit{Connes-Kreimer Hopf algebra.} 
    The coproduct $\delta_{CK}$ of $\kappa[\Lfr^*]$ can be described using the notion of \textit{admissible cut of $t\in \mathcal{T}$}, i.e. a subset (possibly empty) $c$ of the set of edges of $t$ such that for any path from the root to any leaf, there is at most one edge of the path contained in $c$. Given an admissible cut $c$ of $t$, we can delete the edges in $c$ to obtain a collection of rooted subtrees. We will call the \textit{trunk of $t$}, denoted by $R_c(t)$, the subtree which contains the root of $t$, and the \textit{pruning of $t$}, denoted by $P_c(t)$, the monomial given by the forest of the remaining rooted subtrees.
\par Hence, the coproduct $\delta_{CK}:\kappa[\Lfr^*] \to \kappa[\Lfr^*]\otimes \kappa[\Lfr^*]$ is defined by the conditions\linebreak $\delta_{CK}(1) = 1\otimes 1$,
\begin{equation}
    \label{eq:CKcop}
    \delta_{CK}(t) = 1\otimes t + \sum_{c \;\mbox{\scriptsize{admissible cut of }} t} R_c(t) \otimes P_c(t),\quad\mbox{for any }\, t\in \Lfr^*,
\end{equation}
and $\delta_{CK}(u_1\cdots u_n) = \delta_{CK}(u_1)\cdots \delta_{CK}(u_n)$, for any monomial $u_1\cdots u_n\in \kappa[\Lfr^*]$.
\\\par It will also be useful to define the following notions associated to non-planar rooted trees.
    
    \begin{defi}
    \begin{enumerate}
        \item Let $t,u_1,\ldots,u_n\in \T$. We define the non-planar rooted tree $B^+(u_1\cdots u_n) $ as the tree obtained by adding a common root to the forest $u_1\cdots u_n$. We also denote $B^-(t)$ the forest $f$ of non-planar rooted trees such that $B^+(f) = t$. The notation is similar to the one used for decorated rooted trees.
        \item We define the \textit{tree factorial of $t\in \T$}, denoted by $t!$, as the integer recursively defined  by $t! = 1$ if $t$ is the single-vertex tree, and if $t = B^+(s_1\cdots s_m)$, then $t! = |t|s_1! \cdots s_m!$. If $f = t_1\cdots t_n$ is a forest of trees $t_1,\ldots,t_n\in \T$, we define $f!: = t_1!\cdots t_n!$.
    \end{enumerate}
    \end{defi}

The duality, however, is slightly more complicated than for the pre-Lie algebra of words since it entails the so-called \textit{internal symmetry factor of a rooted tree $t$}, defined as $\sigma(t) := |\operatorname{Aut}(t)\,|$. It is obtained as follows \cite[Thm. 301 A]{Butcher}. Write $t=B^+(t_1^{m_1}\cdots t_p^{m_p})$ if the tree $t$ is obtained by grafting on its root the trees $t_1,\ldots,t_p$, all distinct and respectively with multiplicities $m_1,\ldots,m_p$. Then,
    set $\sigma(\bullet):=1$ and define inductively $\sigma(t)$ by
    $$\sigma(t):=\prod\limits_{i=1}^p m_i!\sigma(t_i)^{m_i}.$$

\begin{theo}[{\cite[Proposition 4.4]{Hofm}}]
The duality between the Hopf algebras $\kappa[\Lfr]$ and $\kappa[\Lfr^*]$ is given by the pairing
\begin{equation}
\label{eq:duality}
    \langle s_1\cdots s_m | t_1\cdots t_n \rangle = \left\{\begin{tabular}{c l}
        $\sigma(B^+(t_1\cdots t_n))$&$ \mbox{if } s_1\cdots s_m \simeq t_1\cdots t_n$,  \\ 
         0 & \mbox{otherwise,}
    \end{tabular} \right.
\end{equation}
for any monomials $s_1\cdots s_m \in \kappa[\Lfr] $ and $t_1\cdots t_n\in \kappa[\Lfr^*]$.
\end{theo}

\begin{rem}
Observe that, in the particular case of a single tree $t\in \Lfr^*$, Equation \eqref{eq:duality} implies that
\begin{equation}
    \label{eq:duality2}
    \langle t | t \rangle = \sigma(B^+(t)) = \sigma(t).
\end{equation}
\end{rem}

We can finally prove the main result of this section. To simplify the notation, if $(L,\lhd)$ is a pre-Lie algebra and $\alpha\in L$, we write
$$\exp^\lhd(\alpha) = \sum_{n\geq 1} \frac{1}{n!} r^{(n)}_{\alpha}(\alpha),$$
where in general, $r^{(0)}_{\gamma}(\alpha) =1$, $r^{(1)}_{\gamma}(\alpha) =\alpha$, and $r^{(n)}_{\gamma}(\alpha) = r^{(n-1)}_{\gamma}(\alpha)\lhd \gamma$, for $n\geq 2$.

\begin{prop}
\label{prop:Exppre-Lie}
For the generator $\bullet$ of the free pre-Lie algebra $\Lfr$, we have that 
\begin{equation}
\label{eq:ExpFreepre-Lie}
    \exp^\lhd(\bullet) = \displaystyle \sum_{t\in \mathcal{T}} \frac{|t|!}{\sigma(t)t!} t.
\end{equation}
\end{prop}

\begin{proof}
Let us expand the pre-Lie exponential in the basis of rooted trees, $\exp^\lhd(\bullet) = \sum\limits_{t\in \T} c(t)t.$ Our objective is to describe the coefficients $c(t)$ for all $t\in \mathcal{T}$. For this purpose, consider a tree $t_i\in \Lfr^*$ with $|t_i|=k$ vertices.  Observe that, by definition of the pre-Lie product, $t_i$ only appears in the expansion of $\frac{r^{(k)}_\bullet (\bullet)}{k!}$. By \eqref{eq:duality2} we have
\begin{eqnarray*}
\langle \exp^\lhd(\bullet) | t_i \rangle &=& \frac{1}{k!}\langle r^{(k)}_\bullet (\bullet)| t_i\rangle
\\ &=& \frac{1}{k!}\langle c(t_i) t_i | t_i\rangle 
\\ &=& \frac{1}{k!} c(t_i)\sigma(t_i).
\end{eqnarray*}
On the other hand, by using Lemma \ref{lem:iterated2} and the forest formula \eqref{eq:forestformula} with the basis $\T = \{t_j\}_{j\in \mathbb{N}^*}$ of $\Lfr^*$,  we obtain that
\begin{eqnarray}
\nonumber\langle r^{(k)}_\bullet (\bullet) | t_i\rangle &=& \langle \bullet \otimes \cdots \otimes \bullet| \delta_{\mathrm{irr}}^{[k]}(t_i)\rangle
\\ \label{eq:aux64} &=& \sum_{\substack{T\in \mathcal{T}_i\\ |T| = k}} \sum_{f\in \operatorname{lin}(T)} \lambda(T) \langle \bullet \otimes \cdots\otimes \bullet |C(f) \rangle,
\end{eqnarray}
where we recall that $\mathcal{T}_i$ is the set of decorated rooted trees associated to $t_i$ (Definition \ref{def:tree}.2). 

By definition of the pairing, the non-vanishing terms correspond {to $T\in \T_i$ such that $\lambda(T)\neq 0$ and}
$$C(f)=\bullet\otimes \dots\otimes\bullet $$
{for any $f\in \operatorname{lin(T)}$.} The latter is the case when $T\in\mathcal{T}_i$ with $d_2(x)=0$ for any vertex $x$ of $T$, {where we recall that $t_0 = \bullet$.} Let us show by induction on the number of vertices of $t_i$ that there is only one tree $T$ {that provides a non-zero contribution in \eqref{eq:aux64}}: the decorated tree denoted $T_i^\bullet$, which is the tree $t_i$ where each vertex $v$
is decorated by the pair $(i_v;0)$, where $t_{i_v}$ is the maximal subtree of $t_i$ whose root is $v$.

We refer freely to the arguments in the proof of Theorem \ref{thm:forestform0}. Write $t_i=B^+(t_{i_1}\dots t_{i_n})$. In the calculation of $\delta_{\mathrm{irr}}^{[k]}(t_i)$ in the proof of Theorem \ref{thm:forestform0}, what we have just observed means that in the first step (the calculation of $\overline{\delta}(t_i)$), the only term to be kept is $\bullet\otimes t_{i_1}\dots t_{i_n}$. The next step involves the calculation of the $\delta^{[k-1]}(t_{i_j})$, to which we can apply the induction hypothesis: the only tree in ${\mathcal T}_{i_j}$ that contributes to the calculation of $\langle r^{(k)}_\bullet (\bullet) | t_i\rangle$ is $T_{{i_j}}^\bullet$. Summing up: the only decorated tree that contributes to the calculation of $\langle r^{(k)}_\bullet (\bullet) | t_i\rangle$ is obtained by grafting the $T_{{i_j}}^\bullet$ to a common root $\bullet$ decorated by $(i;0)$: this is precisely $T_i^\bullet$, as expected. 

To resume, the only relevant decorated tree  is $T_i^\bullet = B^+_{(i;0)}(T_{{i_1}}^\bullet\cdots T_{{i_n}}^\bullet)$. We can rewrite the forest $F=T_{{i_1}}^\bullet\cdots T_{{i_n}}^\bullet$ as in \eqref{eq:multiforest} as
$$ \{(T_{a_1}^{\bullet})^{k_1}\}\cup\cdots\cup \{ (T_{a_p}^{\bullet})^{k_p}\}$$
where the $a_i$ are distinct indexes. This implies $\operatorname{sym}(F) = 1$. {Also, it is clear that $\lambda^{i;0}_{i_1,\ldots,i_n}=1$.} By induction on the subtrees of $T_i^\bullet$, to which the same argument applies, we conclude that
{$$\lambda(T_i^\bullet) = \lambda^{i;0}_{i_1,\ldots,i_n}\operatorname{sym}(F)\lambda(T_{i_1}^\bullet)\cdots \lambda(T_{i_n}^\bullet)=1.$$} Therefore

\begin{eqnarray*}
\sum_{\substack{T\in \mathcal{T}_i\\ |T| = k}} \sum_{f\in \operatorname{lin}(T)} \lambda(T) \langle \bullet \otimes \cdots\otimes \bullet |C(f) \rangle &=&  \sum_{f\in \operatorname{lin}(T_i^\bullet)} \langle \bullet \otimes \cdots\otimes \bullet |C(f) \rangle
\\ &=&  \sum_{f\in \operatorname{lin}(t_i)} \prod_{i=1}^k \langle \bullet |\bullet\rangle
\\ &=& \sum_{f\in \operatorname{lin}(t_i)}1
\\ &=& m(t_i),
\end{eqnarray*}
where $m(t_i)= |\operatorname{lin}(t_i)|$. Finally, we can compare it with the expression obtained at the beginning to conclude that
$$c(t)\sigma(t) = m (t).$$ 
By using the well-known fact $m(t) = \frac{|t|!}{t!},$ (\cite[Proposition 3.3]{AHLV}), it follows that
$$c(t) = \frac{m(t)}{\sigma(t)} = \frac{|t|!}{\sigma(t)t!},$$
as we wanted to show.
\end{proof}

\begin{rem}
As was stated at the beginning of the present section, the formula pre-Lie exponential of the generator of the free pre-Lie algebra is a classical result. The coefficients 
$$\mathrm{CM}(t) := \frac{|t|!}{\sigma(t)t!}$$
for $t\in \mathcal{T}$ appear in various contexts and are known as the \textit{Connes-Moscovici coefficients}; there exists a direct method to compute them  (\cite{Brouder}).
\end{rem}

\section{Murua's Coefficients and the Magnus Operator I}
\label{sec:MuruaMagnusI}
In this section and the following one, $(L,\lhd)$ will denote a connected graded and locally finite dimensional pre-Lie algebra.

The Magnus operator $\Omega$ is the compositional inverse of $\exp^{\lhd}$: $$\exp^{\lhd}(\Omega(\alpha)) = \alpha = \Omega(\exp^\lhd(\alpha)),$$
for any $\alpha \in L$. See e.g. \cite[Chap. 6]{CP21}. When used in the context of matrix-valued linear differential equations, the Magnus operator leads to the Magnus formula that computes the derivative of the logarithm of the fundamental solution of such an equation \cite[Remark 6.5.3]{CP21}. Recall that the Baker-Campbell-Hausdorff formula computes the logarithm of the fundamental solution of such an equation. 

\par The coefficients associated to a rooted tree defined in this section have appeared in Murua's analysis of the continuous Baker-Campbell-Hausdorff problem in a Hall basis  (\cite{Mu06}), and, more recently, as the coefficients defining the cumulant-moment and cumulant-cumulant formulas associated to monotone independence in the context of non-commutative probability theory (\cite{cefpp21}).
We feature once again that one of our aims was to deepen the understanding of Murua's article and its relations to our work in non-commutative probability: the methods we develop hereafter shed a new light on his formulas and constructions on (non-planar) trees by showing that they emerge naturally from our forest formulas. Let us mention that there exist other approaches to the Magnus formula by means of (planar) trees: see \cite{iserles} and \cite{efm}.

Recall first some general results and notations. An \textit{ordered partition of $[n]$ of length $k$} is a $k$-tuple $(I_1,\ldots,I_k)$ where the $I_1,\ldots,I_k$ are pairwise disjoint non-empty subsets of $[n]$ such that $I_1\cup\cdots \cup I_k = [n]$. The set of ordered partitions of $[n]$ of length $k$ will be denoted by $\mathcal{OP}^k(n)$.
\par Now, if $L$ is connected graded and locally finite dimensional pre-Lie algebra, consider $b_1,\ldots,b_n \in L$. Denote $\hat{L}$ its completion with respect to the grading. For the monomial $b_1\cdots b_n\in \kappa[L]$, we define 
\begin{equation}
    \operatorname{sol}_1(b_1\cdots b_n) := \sum_{k=1}^n \frac{(-1)^{k-1}}{k} \sum_{\substack{\pi \in \mathcal{OP}^k(n)\\ \pi = (I_1,\ldots,I_k)}} b_{I_1} * \cdots * b_{I_k},
\end{equation}
where $*$ stands for the associative product \eqref{eq:product} in the Hopf algebra $(\kappa[L],*,\Delta)$. The morphism $\operatorname{sol}_1$ is extended linearly to $L$ and to its completion $\hat{L}$.

\begin{theo}[{\cite[Theorem 4.2]{CHP13}}]
\label{thm:AltMagnus}
For $a\in L$, we have, in $\hat{L}$:
\begin{equation}
    \Omega(a) = \operatorname{sol}_1(\exp(a)),
\end{equation}
where $\exp(a)=\sum_{n\geq0} \frac{a^{\cdot n}}{n!}$, with $a^{\cdot n}$ the $n$-th usual polynomial power of $a$ in $\kappa[{L}]$ (we insist that it is {\rm not} defined using the $\ast$ product).
\end{theo}

\begin{defi}
\label{def:murua}
Let $t$ be a non-planar rooted tree with  $|t|$ vertices. For any integer $0 < k \leq |t|$, we denote by $\omega_k(t) := |k-\operatorname{lin}(t)|$ the number of surjective, strictly order preserving maps $f : t \to [k]$. The \textit{Murua's coefficient of $t$} is defined by
\begin{equation}
    \omega(t) := \sum_{k=1}^{|t|} \frac{(-1)^{k-1}}{k}\omega_k(t).
\end{equation}
If $f = t_1\cdots t_s$ is a forest with trees $t_1,\ldots,t_s$, we set $\omega(f) : = \omega(t_1)\cdots \omega(t_s).$
\end{defi}
Notice that, besides in the work of Murua, the coefficients and the following formula had also appeared (in a different framework) in the work of P. Chartier, E. Hairer and G. Vilmart, namely in the context of the numerical analysis of PDEs \cite[Remark 12]{Mu06}. We connect here these works with our forest formulas for coproducts and feature that our approach is not limited to the study of the free pre-Lie algebra.

\begin{prop}
\label{prop:MagnusFree}
For the generator $\bullet$ of the free pre-Lie algebra $\Lfr$, we have that
\begin{equation}
    \label{eq:MagnusFreepre-Lie}
    \Omega(\bullet) = \sum_{t\in \mathcal{T}} \frac{\omega(t)}{\sigma(t)} t.
\end{equation}
\end{prop}

\begin{proof}
The pre-Lie Magnus operator action on the generator of the free pre-Lie algebra can be written $\Omega:= \Omega(\bullet) = \sum_{t\in \mathcal{T}}d(t)t$, for some coefficients $d(t)\in \kappa$. Let $t_i\in \T$. We will show that $d(t_i) = \omega(t_i)/\sigma(t_i).$ On one hand, \eqref{eq:duality2} implies that
$$\langle \Omega|t_i \rangle = d(t_i)\sigma(t_i).$$
On the other hand, by Theorem \ref{thm:AltMagnus} we also have that
\begin{equation}
    \label{eq:auxnew}
    \langle \operatorname{sol}_1(\exp(\bullet)) |t_i\rangle = d(t_i)\sigma(t_i).
\end{equation}
We will analyse the left-hand side of the previous equation. We know that $\exp(\bullet) = \sum_{n\geq0} \frac{\bullet ^n}{n!}$, where $\bullet^n$ stands for the monomial given by the polynomial product of $n$ single-vertex trees. Thus, we can write
\begin{equation}
\label{eq:solbullet}
    \sol(\bullet^n) = \sum_{k=1}^n \frac{(-1)^{k-1}}{k} \sum_{\substack{j_1,j_2\ldots,j_k\geq1\\ j_1+j_2+\cdots +j_k = n}} {n\choose j_1,j_2,\ldots,j_k} \bullet^{j_1} * \bullet^{j_2} *\cdots *\bullet^{j_k}.
\end{equation}
By the definition of the product $*$, the term $\bullet^{j_1}*\cdots * \bullet^{j_k}$ in the equation above is a forest with $j_1+\cdots + j_k = n$ vertices. Hence, if $|t_i| = n$ we have that
\begin{equation}
    \label{eq:auxnew2}
    \langle \sol(\exp(\bullet))|t_i \rangle = \frac{1}{n!} \langle \sol(\bullet^n)|t_i\rangle.
\end{equation} 

Furthermore, it is easy to get from \eqref{eq:product} that when $j_1\geq2$, then $\bullet^{j_1}*\left( \bullet^{j_2}*\cdots * \bullet^{j_k}\right)$ is a linear combination of forests with more than one tree, and this will produce a zero contribution in $\langle \sol(\bullet^n)|t_i\rangle.$ Moreover, we have that the only term in the right-hand side of \eqref{eq:product} for the product $\bullet* \bullet^{j_2}$, that produces a non-zero contribution in \eqref{eq:auxnew2}, is given when $B_0 = 1$, and this term is precisely  $$ \bullet\{\underbrace{\bullet,\ldots,\bullet}_{j_2\,\mbox{\scriptsize{ times}}}\}=: \bullet\{\bullet^{j_2}\}\in L_\mathrm{free},$$ where the braces are associated to the pre-Lie product \eqref{eq:pre-LieTrees}. Since $\bullet\{\bullet^{j_2}\}$ is again an element in $L_\mathrm{free}$, we will have that the only term in $ (\cdots((\bullet * \bullet^{j_2})* \bullet^{j_3})*\cdots ) *\bullet^{j_k}$ that produces a non-zero contribution in \eqref{eq:auxnew2} is precisely
$$(\cdots((\bullet\{\bullet^{j_2}\})\{\bullet^{j_3}\})\cdots)\{\bullet^{j_k}\}. $$
Thus, from \eqref{eq:solbullet} and Lemma \ref{lem:iteratedbrace}, we can compute the action of the iterated brace products and obtain
\begin{eqnarray}
\nonumber
 \frac{1}{n!} \langle \sol(\bullet^n)|t_i\rangle &=& \frac{1}{n!} \sum_{k=1}^n \frac{(-1)^{k-1}}{k} \sum_{\substack{j_2\ldots,j_k\geq1\\ 1+j_2+\cdots +j_k = n}} {n\choose 1,j_2,\ldots,j_k} \langle (\cdots(\bullet\{\bullet^{j_2}\})\cdots)\{\bullet^{j_k}\} |  t_i \rangle 
 \\ &=& \sum_{k=1}^n\frac{(-1)^{k-1}}{k} \sum_{\substack{j_2\ldots,j_k\geq1\\ 1+j_2+\cdots +j_k = n}} \frac{1}{j_2!\cdots j_k!} \langle  \bullet \otimes \bullet^{j_2}\otimes\cdots \otimes \bullet^{j_k} | \overline{\delta}^{[k]}(t_i)\rangle.
 \label{eq:auxsol2}
\end{eqnarray}
Now, by using the forest formula \eqref{eq:forestformula0}, it follows that
\begin{eqnarray*}
\langle  \bullet \otimes \bullet^{j_2}\otimes\cdots \otimes \bullet^{j_k} | \overline{\delta}^{[k]}(t_i)\rangle &=& \sum_{T\in \T_i}\sum_{f\in k-\operatorname{lin}(T)} \lambda(T)\langle \bullet \otimes \bullet^{j_2}\otimes\cdots \otimes \bullet^{j_k} |C(f)\rangle.
\end{eqnarray*}
Observe that the decorated trees $T\in \T_i$ that produce a non-zero contribution must satisfy that $|T| = 1+j_2+\cdots+j_k = n$ and $d_2(x) = 0$ for any $x\in V(T)$. Thus, as in the proof of Proposition \ref{prop:Exppre-Lie}, the only tree $T\in \T_i$ that may produce a non-zero contribution in the above equation is $T^\bullet_i = t_i$ with $d(x) = (i_x;0)$ for any $x\in V(T_i^\bullet)$, where $t_{i_x}$ is the maximal subtree of $t_i$ whose root is $x$. 
\par Recall that we also have that $\lambda(T_i^\bullet) = 1$. We get
$$
\sum_{\substack{j_2\ldots,j_k\geq1\\ 1+j_2+\cdots +j_k = n}} \frac{1}{j_2!\cdots j_k!} \langle  \bullet \otimes \bullet^{j_2}\otimes\cdots \otimes \bullet^{j_k} | \overline{\delta}^{[k]}(t_i)\rangle $$ $$ = \sum_{f\in k-\operatorname{lin}(T_i^\bullet)} \sum_{\substack{j_2\ldots,j_k\geq1\\ 1+j_2+\cdots +j_k = n}} \frac{1}{j_2!\cdots j_k!} \langle  \bullet \otimes \bullet^{j_2}\otimes\cdots \otimes \bullet^{j_k} | C(f) \rangle.
$$
Notice that, given $f\in k-\operatorname{lin}(T_i^\bullet) = k-\operatorname{lin}(t_i)$, there is exactly one tuple $(j_2,\ldots,j_k)$ such that $\langle  \bullet \otimes \bullet^{j_2}\otimes\cdots \otimes \bullet^{j_k} | C(f) \rangle \neq 0$  given by $j_m = |f^{-1}(m)|=:j_m'$, for any $2\leq m\leq k$. Hence, the right-hand side of the above equation is equal to
\begin{eqnarray*}
\sum_{f\in k-\operatorname{lin}(t_i)} \frac{1}{j_2'!\cdots j_k'!} \langle  \bullet \otimes \bullet^{j_2'}\otimes\cdots \otimes \bullet^{j_k'} | C(f) \rangle &=& \sum_{f\in k-\operatorname{lin}(t_i)}  \frac{1}{j_2'!\cdots j_k'!} \prod_{m=2}^k \langle \bullet^{j_m'}| \bullet^{j_m'}\rangle
\\ &=& \sum_{f\in k-\operatorname{lin}(t_i)} \frac{1}{j_2'!\cdots j_k'!} \prod_{m=2}^k \sigma(B^+(\bullet^{j_m'})) 
\\ &=& \sum_{f\in k-\operatorname{lin}(t_i)} \frac{1}{j_2'!\cdots j_k'!} \prod_{m=2}^k j_m'!
\\&=& \sum_{f\in k-\operatorname{lin}(t_i)} 1
\\ &=& |k-\operatorname{lin}(t_i)|,
\end{eqnarray*}
where in the second equality, we used \eqref{eq:duality}. Finally, recalling that $\omega_k(t_i) = |k-\operatorname{lin}(t_i)|$ as in Definition \ref{def:murua}, combining the previous development with \eqref{eq:auxsol2}, we obtain that
\begin{eqnarray*}
\frac{1}{n!}\langle\sol(\bullet^n)|t_i\rangle &=& \sum_{k=1}^n \frac{(-1)^{k-1}}{k} \omega_k(t_i) = \omega(t_i). 
\end{eqnarray*}
Therefore, using \eqref{eq:auxnew} we conclude that $$d(t_i) = \frac{\omega(t_i)}{\sigma(t_i)},$$
as we wanted to show.
\end{proof}

\section{Murua's Coefficients and the Magnus Operator II}
\label{sec:MuruaMagnusII}
The approach we have developed to the Magnus operator in the previous section is, at out best knowledge, new in that it relies, besides on the forest formula for iterated coproducts, on the particular form taken by the canonical projection from the enveloping algebra of a pre-Lie algebra to the underlying pre-Lie algebra.

In this section, we connect forest formulas for iterated {\it reduced} coproducts with a more standard approach to the Magnus operator, namely the one by means of a recursion  --that is, the Agrachev-Gamkrelidze fixed-point equation for the Magnus operator, see Proposition \ref{eq:MagnusExpansion} below. The resulting proofs and arguments are more involved than the ones in the previous section. We thought it important to include them as they connect explicitly to the arguments in \cite{Mu06}, but also because of the practical use that can be made of this fixed-point equation in numerical analysis.

The following proposition allows to give a recursive and alternative definition of Murua's coefficients. It is in \cite{Mu06} but will also follow from our developments in this section of the article. 

\begin{prop}[{\cite[Remark 11]{Mu06}}]
\label{prop:recursiveMurua}
For any  non-planar rooted tree $t$ with $|t|>1$, we have that
\begin{equation}
    \label{eq:MuruaRecursion}
    \omega(t) = \sum_{s \in K(B^-(t))} \frac{B_{|s|}}{s!} \omega(C^s(B^-(t))),
\end{equation}
where, for $f$ a forest, $K(f)$ stands for the multiset of subforests of $f$ that contain all the roots of the trees of $f$, and where, if $s\in K(f)$, $C^s(f)$ stands for the forest obtained from $f$ by removing the edges that connect the vertices of $s$ with their parents.
\end{prop}
In the previous proposition, forests $f$ and trees $t$ are viewed as posets (with roots as minimal elements). A subforest of $f$ or $t$ is a subset of $f$ or $t$ equipped with the induced order. As usual, subtrees are subforest with a unique minimal element (the root). Notice that $C^s(f)$ has the same set of vertices as $f$ (in the definition of $C^s(f)$, edges are removed, but not vertices).

\begin{prop}
\label{eq:MagnusExpansion}
 The \textit{pre-Lie Magnus operator} $\Omega:L\to L$ is given by
$$\Omega(\alpha) = \sum_{n\geq 0 } \frac{B_n}{n!} r^{(n)}_{\Omega(\alpha)}(\alpha),$$
for any $\alpha\in L$ and $\{B_n\}_{n\geq0}$ is the sequence of Bernoulli numbers.
\end{prop}

Our analysis will rely on the existence of a relation between the families $K(f)$, $C^s(f)$, decorated trees and iterated coproducts explained in the next lemma.

\begin{lemma}
\label{lem:surjectiveMurua}
Let $t_i\in \T$ with $i\neq0$ and denote by $\T_i^r$ to the subset of decorated trees $T\in \T_i$ such that $\lambda(T)\neq0$ and $d_2(\operatorname{root}(T)) = 0$. There is then a surjective map $B:K(B^-(t_i))\to \T_i^r$ such that $|\{s\in K(B^-(t_i))\,:\, B(s)=T\}| = \lambda(T)$, for any $T\in \T_i^r$.
\end{lemma}
\begin{proof}
First of all, let us describe the map $B$. To this end, assume that $t_i = B^+(t_{i_1}\cdots t_{i_k})$ and consider an element $s\in K(B^-(t_i))$ such that $s = t_{j_1}\cdots t_{j_k}$, where for each $1\leq m\leq k$, $t_{j_m}$ is a subtree of $t_{i_m}$, with the same root. Then, we define the decorated tree $B(s) := B^+_{(i;0)}(t_{j_1}\cdots t_{j_k})$, where, for each $1\leq m\leq k$, the decoration of $x\in V(t_{j_m})$ is given by:
\begin{itemize}
    \item $d_1(x)$ is the index associated to the subtree of $t$ defined by $x$ and all its descendants,
    \item $d_2(x)$ is the index associated to the tree obtained from $t_{d_1(x)}$ once all the subtrees determined by the elements of $V(t_{j_m})\backslash\{x\}$ and their descendants have been removed. Observe that we also delete the edges connecting these elements with their parents in $t_{d_1(x)}$.
\end{itemize}
This map is well-defined since, by construction, $B(s)\in \T_i^r$. 
\par With the purpose of showing that $B$ is a surjective map, we take a decorated tree $T \in \T_i^r$. By the proof of Proposition \ref{prop:Exppre-Lie}, we can write $T = B^+_{(i;0)}(T_1\cdots T_k)$, where  for $1\leq m\leq k$, $T_m$ is a decorated tree associated to $t_{i_m}$.  Also, the condition $\lambda(T)\neq0$ implies that $\lambda^{d(x)}_{d_1(\operatorname{succ}(x))}\neq 0$ for any $x\in V(T)$. Recall that by the definition of the Connes-Kreimer coproduct, $\lambda^{d(x)}_{d_1(\operatorname{succ}(x))}$ is the number of admissible cuts $c$ of $t_{d_1(x)}$ such that $R_c(t_{d_1(x)}) = t_{d_2(x)}$ and $P_c(t_{d_1(x)}) = t_{d_1(y_1)}\cdots t_{d_1(y_r)}$ with $\operatorname{succ}(x) = \{y_1,\ldots,y_r\}$. 
\par
Now, observe that, for $1\leq m\leq k$, the process of construction of $T_m$ from $t_{i_m}$ corresponds to a certain contraction of $t_{i_m}$. To be more precise, every $x\in V(T_m)$ can be seen as a vertex in $t_{i_m}$. Then, $x\in V(T_m)$ is obtained by collapsing to a single vertex the subtree $t_{d_2(x)}$ in $t_{i_m}$ which has $x$ as root. We denote $\lambda'(T_m)$ the number of ways in which we can obtain $T_m$ from $t_{i_m}$ from the above procedure. Since $V(T_m)$ must contain the root of $t_{i_m}$, it is clear that if $|T_m| = 1$, then $T_m$ is the single-vertex tree associated to $t_{i_m}$. Thus we have that $\lambda'(T_m) = 1$ in this case.
\par Now assume that $|T_m|>1$, so $B^-(T_m) = T_{m_1}\cdots T_{m_\ell}$ is a non-empty forest, where $T_{m_j}$ is associated to a $t_{i_{m_j}}$, for each $1\leq j\leq \ell$. 
Observe that the contraction process can be equivalently considered as making admissible cuts. More precisely, if $x$ stands for the root of $T_m$, there are $\lambda^{d(x)}_{i_{m_1},\ldots,i_{m_\ell}}$ ways to obtain $t_{d_2(x)}$ as a subtree of $t_{i_m}$, with $x$ as root, in a such a way that when we collapse $t_{d_2(x)}$ to a vertex, the remaining children of $x$ are the roots of $t_{i_{m_1}},\ldots,t_{i_{m_k}}$. In addition, a symmetry factor appears: there are $\operatorname{sym}(B^-(T_m))$ ways to allocate the decorated trees of the decorated forest $T_{m_1}\cdots T_{m_\ell}$ to the subtrees $t_{i_{m_j}}$. Since the roots of $T_{m_1},\ldots,T_{m_\ell}$ are elements in $V(T_m)$, we can proceed recursively to deduce that the number of ways of obtaining $T_m$ from $t_{i_m}$ by the contraction process is given by
$$\lambda'(T_m) = \lambda^{d(x)}_{i_{m_1},\ldots,i_{m_\ell}}\operatorname{sym}(B^{-}(T_m))\lambda'(T_{m_1})\cdots \lambda'(T_{m_\ell}).$$
Hence, we get that $\lambda'(T_m) = \lambda(T_m)$, for any $1\leq m\leq k$. From this we conclude that the number $\lambda(T) = \lambda^{i,0}_{i_1,\ldots,i_k}\operatorname{sym}(B^-(T))\lambda(T_1)\cdots\lambda(T_k)$ counts the number of ways in which we can contract the subtrees $t_{i_1},\ldots,t_{i_k}$ in order to obtain the trees $T_1,\ldots,T_k$. In other words, there are $\lambda(T)$ elements  $s\in K(B^{-}(t_i))$ such that $B(s) = T$, as we wanted to prove.
\end{proof}

\begin{ex}
Consider the tree 
$$t_i = 
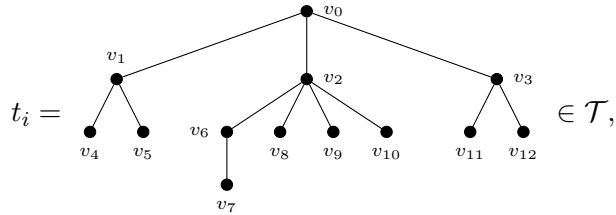
\begin{tikzpicture}[baseline={([yshift=-1ex]current bounding box.center)},scale=1,
level 1/.style={level distance=9mm,sibling distance=25mm},
level 2/.style={level distance=7mm,sibling distance=7mm},
level 3/.style={level distance=7mm,sibling distance=7mm}]
\node(0)[solid node,label=right:{\tiny $v_0$}]{} 
child{node(1)[solid node,label={[font=\tiny,text=black]:$v_1$}]{}
    child{[black] node(11)[solid node,label=below:{\tiny $v_4$}]{}}
	child{[black] node(12)[solid node,label=below:{\tiny $v_5$}]{}}
}
child{node(2)[solid node,label=right:{\tiny $v_2$}]{}
	child{[black] node(21)[solid node,label=left:{\tiny $v_6$}]{}
	    child{[black] node(211)[solid node,label=below:{\tiny $v_7$}]{}}
	}
	child{[black] node(22)[solid node,label=below:{\tiny $v_8$}]{}}
	child{[black] node(23)[solid node,label=below:{\tiny $v_9$}]{}}
	child{[black] node(24)[solid node,label=below:{\tiny $v_{10}$}]{}}
}
child{node(3)[solid node,label=right:{\tiny $v_3$}]{}
	child{[black] node(31)[solid node,label=below:{\tiny $v_{11}$}]{}}
	child{[black] node(32)[solid node,label=below:{\tiny $v_{12}$} ]{}}
};
\end{tikzpicture} \in \T,$$ 
and consider the decorated tree
$$ T = 
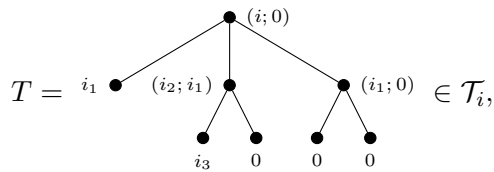
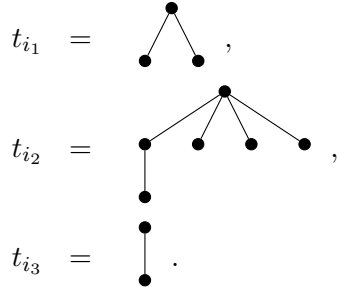
\begin{tikzpicture}[baseline={([yshift=-1ex]current bounding box.center)},scale=1,
level 1/.style={level distance=9mm,sibling distance=15mm},
level 2/.style={level distance=7mm,sibling distance=7mm},
level 3/.style={level distance=7mm,sibling distance=7mm}]
\node(0)[solid node,label=right:{\tiny $(i;0)$}]{} 
child{node(1)[solid node,label=left:{\tiny $i_1$}]{}
}
child{node(2)[solid node,label=left:{\tiny $(i_2;i_1)$}]{}
	child{[black] node(21)[solid node,label=below:{\tiny $i_3$}]{}
	}
	child{[black] node(22)[solid node,label=below:{\tiny $0$}]{}}
}
child{node(3)[solid node,label=right:{\tiny $(i_1;0)$}]{}
	child{[black] node(31)[solid node,label=below:{\tiny $0$}]{}}
	child{[black] node(32)[solid node,label=below:{\tiny $0$} ]{}}
};
\end{tikzpicture} \in \T_i,$$
where
\begin{eqnarray*}
t_{i_1} &=&   \begin{tikzpicture}[baseline={([yshift=-1ex]current bounding box.center)},scale=1,
level 1/.style={level distance=7mm,sibling distance=7mm}]
\node(0)[solid node,label=right:{\tiny }]{} 
child{node(1)[solid node,label=left:{\tiny }]{}
}
child{node(2)[solid node,label=right:{\tiny }]{}
};
\end{tikzpicture},
\\ t_{i_2} &=&\begin{tikzpicture}[baseline={([yshift=-1ex]current bounding box.center)},scale=1,
level 1/.style={level distance=7mm,sibling distance=7mm}]
\node(0)[solid node,label=right:{\tiny }]{} 
child{node(1)[solid node,label=left:{\tiny }]{}
    child{[black] node(11)[solid node,label=right:{\tiny }]{}}
}
child{node(2)[solid node,label=right:{\tiny }]{}}
child{node(3)[solid node,label=right:{\tiny }]{}}
child{node(4)[solid node,label=right:{\tiny }]{}
};
\end{tikzpicture},
\\ t_{i_3} &=&\begin{tikzpicture}[baseline={([yshift=-1ex]current bounding box.center)},scale=1,
level 1/.style={level distance=7mm,sibling distance=7mm}]
\node(0)[solid node,label=right:{\tiny }]{} 
child{node(1)[solid node,label=left:{\tiny }]{}
};
\end{tikzpicture}.
\end{eqnarray*}
Observe that we have labeled the vertices in $t_i$ since the definition of $K(B^-(t_i))$ requires to consider subposets of $t_i$ --and not just isomorphism classes of subposets. 
\par One can easily compute that $\lambda(T) = \operatorname{sym}(B^-(T)) \lambda^{i_2;i_1}_{i_3,0} = 2\cdot 3 = 6$. The 6 elements in $K(B^-(t_i))$ are depicted in red as subposets of $t_i$ as follows:
$$
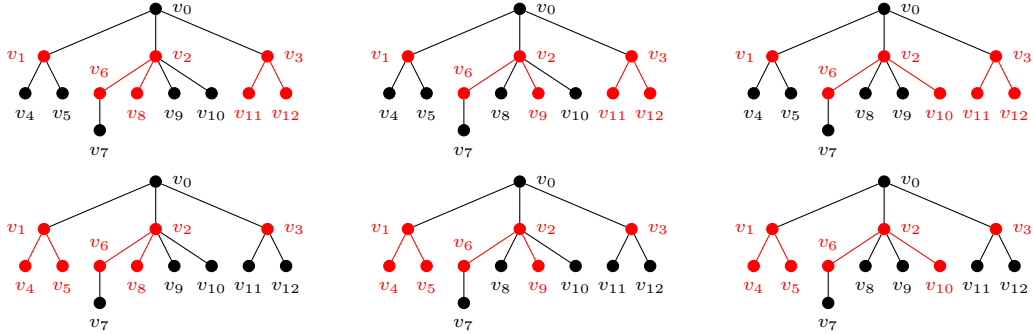
\begin{tabular}{c c c}
   $$ \begin{tikzpicture}[baseline={([yshift=-1ex]current bounding box.center)},scale=0.7,
level 1/.style={level distance=9mm,sibling distance=21mm},
level 2/.style={level distance=7mm,sibling distance=7mm},
level 3/.style={level distance=7mm,sibling distance=7mm}]
\node(0)[solid node,label=right:{\tiny $v_0$}]{} 
child{node(1)[solid node,color=red,label=left:{\tiny $\color{red}v_1$}]{}
    child{[black] node(11)[solid node,label=below:{\tiny $v_4$}]{}}
	child{[black] node(12)[solid node,label=below:{\tiny $v_5$}]{}}
}
child{node(2)[solid node,color=red,label=right:{\tiny $\color{red}v_2$}]{}
	child{[red] node(21)[solid node,color=red,label=above:{\tiny $v_6$}]{}
	    child{[black] node(211)[solid node,label=below:{\tiny $v_7$}]{}}
	}
	child{[red] node(22)[solid node,color=red,label=below:{\tiny $v_8$}]{}}
	child{[black] node(23)[solid node,label=below:{\tiny $v_9$}]{}}
	child{[black] node(24)[solid node,label=below:{\tiny $v_{10}$}]{}}
}
child{node(3)[solid node,color=red,label=right:{\tiny $\color{red}v_3$}]{}
	child{[red] node(31)[solid node,color=red,label=below:{\tiny $v_{11}$}]{}}
	child{[red] node(32)[solid node,color=red,label=below:{\tiny $v_{12}$} ]{}}
};
\end{tikzpicture} $$  & $$ \begin{tikzpicture}[baseline={([yshift=-1ex]current bounding box.center)},scale=0.7,
level 1/.style={level distance=9mm,sibling distance=21mm},
level 2/.style={level distance=7mm,sibling distance=7mm},
level 3/.style={level distance=7mm,sibling distance=7mm}]
\node(0)[solid node,label=right:{\tiny $v_0$}]{} 
child{node(1)[solid node,color=red,label=left:{\tiny $\color{red}v_1$}]{}
    child{[black] node(11)[solid node,label=below:{\tiny $v_4$}]{}}
	child{[black] node(12)[solid node,label=below:{\tiny $v_5$}]{}}
}
child{node(2)[solid node,color=red,label=right:{\tiny $\color{red}v_2$}]{}
	child{[red] node(21)[solid node,color=red,label=above:{\tiny $v_6$}]{}
	    child{[black] node(211)[solid node,label=below:{\tiny $v_7$}]{}}
	}
	child{[black] node(22)[solid node, label=below:{\tiny $v_8$}]{}}
	child{[red] node(23)[solid node, color=red,label=below:{\tiny $v_9$}]{}}
	child{[black] node(24)[solid node,label=below:{\tiny $v_{10}$}]{}}
}
child{node(3)[solid node,color=red,label=right:{\tiny $\color{red}v_3$}]{}
	child{[red] node(31)[solid node,color=red,label=below:{\tiny $v_{11}$}]{}}
	child{[red] node(32)[solid node,color=red,label=below:{\tiny $v_{12}$} ]{}}
};
\end{tikzpicture} $$ & $$ \begin{tikzpicture}[baseline={([yshift=-1ex]current bounding box.center)},scale=0.7,
level 1/.style={level distance=9mm,sibling distance=21mm},
level 2/.style={level distance=7mm,sibling distance=7mm},
level 3/.style={level distance=7mm,sibling distance=7mm}]
\node(0)[solid node,label=right:{\tiny $v_0$}]{} 
child{node(1)[solid node,color=red,label=left:{\tiny $\color{red}v_1$}]{}
    child{[black] node(11)[solid node,label=below:{\tiny $v_4$}]{}}
	child{[black] node(12)[solid node,label=below:{\tiny $v_5$}]{}}
}
child{node(2)[solid node,color=red,label=right:{\tiny $\color{red}v_2$}]{}
	child{[red] node(21)[solid node,color=red,label=above:{\tiny $v_6$}]{}
	    child{[black] node(211)[solid node,label=below:{\tiny $v_7$}]{}}
	}
	child{[black] node(22)[solid node,label=below:{\tiny $v_8$}]{}}
	child{[black] node(23)[solid node,label=below:{\tiny $v_9$}]{}}
	child{[red] node(24)[solid node,color=red,label=below:{\tiny $v_{10}$}]{}}
}
child{node(3)[solid node,color=red,label=right:{\tiny $\color{red}v_3$}]{}
	child{[red] node(31)[solid node,color=red,label=below:{\tiny $v_{11}$}]{}}
	child{[red] node(32)[solid node,color=red,label=below:{\tiny $v_{12}$} ]{}}
};
\end{tikzpicture} $$\\
     $$ \begin{tikzpicture}[baseline={([yshift=-1ex]current bounding box.center)},scale=0.7,
level 1/.style={level distance=9mm,sibling distance=21mm},
level 2/.style={level distance=7mm,sibling distance=7mm},
level 3/.style={level distance=7mm,sibling distance=7mm}]
\node(0)[solid node,label=right:{\tiny $v_0$}]{} 
child{node(1)[solid node,color=red,label=left:{\tiny $\color{red}v_1$}]{}
    child{[red] node(11)[solid node,color=red,label=below:{\tiny $v_4$}]{}}
	child{[red] node(12)[solid node,color=red,label=below:{\tiny $v_5$}]{}}
}
child{node(2)[solid node,color=red,label=right:{\tiny $\color{red}v_2$}]{}
	child{[red] node(21)[solid node,color=red,label=above:{\tiny $v_6$}]{}
	    child{[black] node(211)[solid node,label=below:{\tiny $v_7$}]{}}
	}
	child{[red] node(22)[solid node,color=red,label=below:{\tiny $v_8$}]{}}
	child{[black] node(23)[solid node,label=below:{\tiny $v_9$}]{}}
	child{[black] node(24)[solid node,label=below:{\tiny $v_{10}$}]{}}
}
child{node(3)[solid node,color=red,label=right:{\tiny $\color{red}v_3$}]{}
	child{[black] node(31)[solid node,label=below:{\tiny $v_{11}$}]{}}
	child{[black] node(32)[solid node,label=below:{\tiny $v_{12}$} ]{}}
};
\end{tikzpicture} $$  & $$ \begin{tikzpicture}[baseline={([yshift=-1ex]current bounding box.center)},scale=0.7,
level 1/.style={level distance=9mm,sibling distance=21mm},
level 2/.style={level distance=7mm,sibling distance=7mm},
level 3/.style={level distance=7mm,sibling distance=7mm}]
\node(0)[solid node,label=right:{\tiny $v_0$}]{} 
child{node(1)[solid node,color=red,label=left:{\tiny $\color{red}v_1$}]{}
    child{[red] node(11)[solid node,color=red,label=below:{\tiny $v_4$}]{}}
	child{[red] node(12)[solid node,color=red,label=below:{\tiny $v_5$}]{}}
}
child{node(2)[solid node,color=red,label=right:{\tiny $\color{red}v_2$}]{}
	child{[red] node(21)[solid node,color=red,label=above:{\tiny $v_6$}]{}
	    child{[black] node(211)[solid node,label=below:{\tiny $v_7$}]{}}
	}
	child{[black] node(22)[solid node, label=below:{\tiny $v_8$}]{}}
	child{[red] node(23)[solid node, color=red,label=below:{\tiny $v_9$}]{}}
	child{[black] node(24)[solid node,label=below:{\tiny $v_{10}$}]{}}
}
child{node(3)[solid node,color=red,label=right:{\tiny $\color{red}v_3$}]{}
	child{[black] node(31)[solid node,label=below:{\tiny $v_{11}$}]{}}
	child{[black] node(32)[solid node,label=below:{\tiny $v_{12}$} ]{}}
};
\end{tikzpicture} $$ & $$ \begin{tikzpicture}[baseline={([yshift=-1ex]current bounding box.center)},scale=0.7,
level 1/.style={level distance=9mm,sibling distance=21mm},
level 2/.style={level distance=7mm,sibling distance=7mm},
level 3/.style={level distance=7mm,sibling distance=7mm}]
\node(0)[solid node,label=right:{\tiny $v_0$}]{} 
child{node(1)[solid node,color=red,label=left:{\tiny $\color{red}v_1$}]{}
    child{[red] node(11)[solid node,color=red,label=below:{\tiny $v_4$}]{}}
	child{[red] node(12)[solid node,color=red,label=below:{\tiny $v_5$}]{}}
}
child{node(2)[solid node,color=red,label=right:{\tiny $\color{red}v_2$}]{}
	child{[red] node(21)[solid node,color=red,label=above:{\tiny $v_6$}]{}
	    child{[black] node(211)[solid node,label=below:{\tiny $v_7$}]{}}
	}
	child{[black] node(22)[solid node,label=below:{\tiny $v_8$}]{}}
	child{[black] node(23)[solid node,label=below:{\tiny $v_9$}]{}}
	child{[red] node(24)[solid node,color=red,label=below:{\tiny $v_{10}$}]{}}
}
child{node(3)[solid node,color=red,label=right:{\tiny $\color{red}v_3$}]{}
	child{[black] node(31)[solid node,color=black,label=below:{\tiny $v_{11}$}]{}}
	child{[black] node(32)[solid node,color=black,label=below:{\tiny $v_{12}$} ]{}}
};
\end{tikzpicture} $$
\end{tabular}
$$
The red-colored edges are the edges cut in the selected admissible cuts to construct out $T$ from $t_i$.
\end{ex}

\begin{rem}
\label{rem:csfMurua}
Let $t_i\in \T$ and $s\in K(B^-(t_i))$. From the definition of $C^s(B^-(t_i))$ and the construction of $B(s)=T$ in the proof of Lemma $\ref{lem:surjectiveMurua}$, one can readily check that \linebreak $C^s(B^-(t_i)) = \{t_{d_2(x)}\,:\, x\in V(T)\backslash \{\operatorname{root}(T)\} \}$.
\end{rem}

After the above technical observations, we can propose another proof of Theorem \ref{prop:MagnusFree}. The reason for including it is, besides the importance of the theorem, the fact that this other proof enlightens new features of the underlying combinatorics of trees and free pre-Lie algebras.

\begin{proof}[Alternative proof of Theorem \ref{prop:MagnusFree}]
Write again $\Omega:=\Omega(\bullet) = \sum\limits_{t\in \mathcal{T}} d(t)t$ for the pre-Lie Magnus operator. We will show that $d(t) = \omega(t)/\sigma(t)$, for every $t\in \T$, using the recursive definitions of the Murua's coefficients and of the Magnus operator. On one hand, \eqref{eq:duality2} implies that 
\begin{equation}
    \label{eq:Mgaux1}
    \langle \Omega |t\rangle= d(t)\sigma(t) =: d'(t).
\end{equation}
for any $t\in \T$. On the other hand, for an element $t\in \T$ we have
\begin{eqnarray*}
\langle \Omega|t\rangle &=& \sum_{m\geq 0} \frac{B_m}{m!} \langle r^{(m)}_{\Omega}(\bullet) | t\rangle
\\ &=& \sum_{m=0}^{|t|-1} \frac{B_m}{m!} \langle \bullet \otimes \Omega\otimes \cdots \otimes \Omega| \delta_{\mathrm{irr}}^{[m+1]}(t)\rangle,
\end{eqnarray*}
where to get the upper limit of the sum in  the last equality, we used that if $s,t$ are trees such that $|s|>|t|$, then $\langle s|t \rangle = 0$.  Notice that if $t=\bullet$, then $\langle \Omega | t\rangle = 1$. If $t\neq \bullet$, by using the forest formula \eqref{eq:forestformula}, the previous equation can be rewritten (with $t=t_i)$ as
\begin{equation}
    \label{eq:Mgaux2}
\langle \Omega|t_i\rangle = \sum_{m=0}^{|t|-1} \frac{B_m}{m!} \sum_{\substack{T\in \T_i\\ |T| = m+1}} \sum_{f\in \operatorname{lin}(T)} \lambda(T) \langle \bullet \otimes \Omega\otimes\cdots \otimes \Omega | C(f) \rangle,
\end{equation}
Observe that, if $T\in \T_i$ is such that $d_2(\operatorname{root}(T)) \neq 0$, then $\langle \bullet \otimes \Omega\otimes \cdots \otimes \Omega|C(f)\rangle = 0$ for any $f\in \operatorname{lin}(T)$. Hence the above sum can be restricted to the set $\T_i^r$ defined in Lemma \ref{lem:surjectiveMurua}, and the right-hand side of \eqref{eq:Mgaux2} can be rearranged as
\begin{equation}
\label{eq:Mgaux3}   
\langle \Omega| t_i\rangle  = \sum_{T\in \T^i_t} \frac{B_{|T|-1}}{(|T|-1)!} m(T)\lambda(T)\prod_{\substack{x\in V(T)\\ x\neq \operatorname{root}(T)}} \langle \Omega | t_{d_2(x)}\rangle.
\end{equation}
Thus, by Lemma \ref{lem:surjectiveMurua}, Remark \ref{rem:csfMurua}, and noticing that for $s\in K(B^-(t_i))$ such that $B(s) =T$ we have $$ m(T) = \frac{|T|!}{T!} = \frac{|T|!}{|T|B^-(T)!} = \frac{(|s|+1)!}{(|s|+1)s!}$$ by definition of tree factorial,  we conclude that the right-hand side of \eqref{eq:Mgaux3} can be expressed as follows:
$$ \sum_{s\in K(B^-(t_i))} \frac{B_{|s|}}{|s|!} \frac{(|s|+1)!}{(|s|+1)s!} \prod_{S\in C^s(B^-(t_i))} d'(S) = \sum_{s\in K(B^-(t_i))} \frac{B_{|s|}}{s!} \prod_{S\in C^s(B^-(t_i))} d'(S).$$
We have shown that $d'(t)$ satisfies the recursion \eqref{eq:MuruaRecursion} with the initial condition $d'(\bullet) = 1$. Therefore $w(t) = d'(t) = \sigma(t)d(t)$, which implies that $d(t) = \omega(t)/\sigma(t)$ for any $t\in \T$, as we wanted to show.
\end{proof}

\section{Cumulant-Cumulant Relations via Forest Formulas}
\label{sec:cumulantformulas}
The objective of this section is to present another application of the forest formula in the context of non-commutative probability, more precisely, in the study of combinatorial relations between the different brands of cumulants for non-commutative independence. 
\par We begin this section by describing the objects used in the combinatorial study of non-commutative cumulants.

\begin{defi}[Non-crossing partitions, irreducible partitions, and the forest of nestings]
\phantom{}
\begin{enumerate}
    \item Let $n\in\mathbb{N}$. A \textit{non-crossing partition of $[n]$} is a collection $\pi = \{V_1,\ldots,V_s\}$ of pairwise disjoint non-empty subsets of $[n]$, called \textit{blocks of $\pi$}, such that $[n] = V_1\cup\cdots\cup V_s$, and, there are no two different blocks $V,W\in \pi$ and $a<b<c<d\in [n]$ such that $a,c\in V$, $b,d\in W$. The set of non-crossing partitions of $[n]$ is denoted by $\NC(n)$. The unique partition with $n$ blocks in $\NC(n)$ is denoted by $0_n$, whilst the unique single-block partition in $\NC(n)$ is denoted by $1_n$.  In general, we define $\NC(X)$ to be the set of non-crossing partitions of a finite totally ordered set $X$, defined in the straightforward way. 
\item We will use several special subsets of non-crossing partitions. An \textit{interval partition of $[n]$} is a partition $\pi\in \NC(n)$ such that every block of $V=\pi$ is of the form $\{i,i+1,\ldots,i+j\}$, for some $1\leq i\leq n$ and $0\leq j\leq n-i$. The set of interval partitions of $[n]$ is denoted by $\operatorname{I}(n)$. Additionally, an \textit{irreducible partition} is a partition $\pi\in \NC(n)$ such that $1$ and $n$ belong to the same block in $\pi$. The set of irreducible partitions of $[n]$ is denoted by $\NCirr(n)$. Moreover, the sets of non-crossing and irreducible partitions with $k$ blocks are denoted $\NC^k(n)$, respectively, $\NCirr^k(n)$.
\item Given $\pi\in \NC(n)$, there is a natural way to endow $\pi$ with a poset structure by defining $V\leq W$, for any two blocks $V,W\in \pi$, if  there are $a,c\in V$ such that $a<b<c$ for any $b\in W$. In other words, $V\leq W$ if $W$ is nested in $V$. If $\pi\in \NCirr(n)$, there is only one block $V\in \pi$ which is minimal: the block that contains $1$ and $n$. In the general case, if  $\pi\in \NC(n)$ and $V_1,\ldots,V_k \in \pi$ are the minimal blocks in $\pi$, we define the \textit{irreducible components of $\pi$} to be the irreducible partitions $\pi_1,\ldots,\pi_k$, where each partition is given by $\pi_j = \{W\in \pi\,:\, V_j\leq W\}$ for $1\leq j\leq k$. It is clear that $\pi$ decomposes uniquely as the union of irreducible partitions $\pi = \pi_1\cup\cdots \cup \pi_k$.
\item Associated to a $\pi \in \NCirr(n)$, we can build a rooted tree $t(\pi)$ with $|\pi|$ vertices. This is done by putting in bijection the vertices of $t(\pi)$ to the blocks of $\pi$ and drawing a directed edge from a vertex $v$ to a vertex $w$ if and only if the corresponding block $V$ associated to $v$ covers the block $W$ associated to $w$, i.e. $V < W$ and there is no other block $U \in \pi$ such that $V < U < W$. The tree obtained by this process is called the \textit{tree of nestings} of $\pi$ and it will be denoted by $t(\pi)$. For the general case $\pi \in \NC(n)$, we can construct the \textit{forest of nestings} of $\pi$ as the ordered forest whose trees are given by the tree of nestings of the irreducible components of $\pi$, where the order of the irreducible components is given by the total order determined by the minimum element of each component.
\end{enumerate}
\end{defi}

\begin{ex} 
 Let $\pi$ and $\sigma$ be non-crossing partitions which are graphically represented by arcs below. The forests of nesting $t(\pi)$ and $t(\sigma)$ are also displayed. In this example, each vertex of the forest is decorated with the minimal element of the associated block in the partition.
\begin{center}
\begin{figure}[H]
$$
\begin{array}{c c c c c c}
\pi &=& \begin{tikzpicture}[baseline={([yshift=0.3ex]current bounding box.center)},thick,font=\small]
     \path 	(0,0) 		node (a) {1}
           	(0.5,0) 	node (b) {2}
           	(1,0) 		node (c) {\phantom{1}}
           	(1.5,0) 	node (d) {4}
           	(2,0) 		node (e) {\phantom{1}}
           	(2.5,0) 	node (f)  {6}
           	(3,0) 		node (g) {\phantom{1}}
           	(3.5,0) 	node (h) {\phantom{1}}
		(4,0) 		node (i) {\phantom{1}}
		(4.5,0) 	node (j) {\phantom{1}};
     \draw (a) -- +(0,0.75) -| (j);
     \draw (b) -- +(0,0.60) -| (b);
     \draw (c) -- +(0,0.75) -| (c);
     \draw (d) -- +(0,0.60) -| (e);
     \draw (e) -- +(0,0.60) -| (h);
     \draw (f) -- +(0,0.50) -| (g);
     \draw (i) -- +(0,0.75) -| (j);
   \end{tikzpicture}, 	& \quad t(\pi) &=&
				  \begin{tikzpicture}[baseline={([yshift=-1ex]current bounding box.center)},scale=0.6,
level 1/.style={level distance=7mm,sibling distance=7mm}]
\node(0)[solid node,label=right:{\tiny 1}]{} 
child{node(1)[solid node,label=left:{\tiny 2}]{}
}
child{node(2)[solid node,label=right:{\tiny 4}]{}
	child{[black] node(11)[solid node,label=right:{\tiny 6}]{}}
};
\end{tikzpicture}
\\[1cm]
\sigma &= & \begin{tikzpicture}[baseline={([yshift=0.3ex]current bounding box.center)},thick,font=\small]
     \path 	(0,0) 		node (a) {1}
           	(0.5,0) 	node (b) {2}
           	(1,0) 		node (c) {\phantom{1}}
           	(1.5,0) 	node (d) {4}
           	(2,0) 		node (e) {\phantom{1}}
           	(2.5,0) 	node (f)  {6}
           	(3,0) 		node (g) {7}
           	(3.5,0) 	node (h) {\phantom{1}};
     \draw (a) -- +(0,0.75) -| (c);
     \draw (b) -- +(0,0.60) -| (b);
     \draw (c) -- +(0,0.75) -| (c);
     \draw (d) -- +(0,0.75) -| (e);
     \draw (e) -- +(0,0.75) -| (h);
     \draw (f)  -- +(0,0.60) -|(f);
     \draw (g) -- +(0,0.60)-|(g); 
   \end{tikzpicture}, 	& \quad t(\sigma) &=&\;\;\begin{tikzpicture}[baseline={([yshift=-1ex]current bounding box.center)},scale=0.6,
level 1/.style={level distance=7mm,sibling distance=7mm}]
\node(0)[solid node,label=right:{\tiny 1}]{} 
child{node(1)[solid node,label=right:{\tiny 2}]{}
};
\end{tikzpicture}
\;\; \begin{tikzpicture}[baseline={([yshift=-1ex]current bounding box.center)},scale=0.6,
level 1/.style={level distance=7mm,sibling distance=7mm}]
\node(0)[solid node,label=right:{\tiny 4}]{} 
child{node(1)[solid node,label=left:{\tiny 6}]{}
}
child{node(2)[solid node,label=right:{\tiny 7}]{}
};
\end{tikzpicture}
\end{array}  
$$
\end{figure}
\end{center}
\end{ex}

In the following, the pair $(\A,\varphi)$ stands for a \textit{non-commutative probability space}, i.e. $(\A,\varphi)$ is a pair such that $\A$ is a unital associative algebra over $\mathbb{C}$ and $\varphi:\A\to\mathbb{C}$ is a linear functional such that $\varphi(1_\A) = 1$, where $1_\A$ is the unit of $\A$. Elements $a\in \A$ are called random variables and $\varphi(a)$ is called the moment of $a$ with respect to $\varphi$.
\par The notions of cumulants for the free, Boolean, and monotone independence can be introduced by the so-called moment-cumulant formulas. Before stating the next definition, we fix the following notation: given a family $\{f_n:\A^n\to\mathbb{C}\}_{n\geq1}$ of multilinear functionals, elements $a_1,\ldots,a_n\in \A$ and $\pi\in \NC(n)$, we write
\begin{equation}
    f_\pi(a_1,\ldots,a_n) := \prod_{V\in \pi} f_{|V|}(a_1,\ldots,a_n|V),
\end{equation}
where $f_{|V|}(a_1,\ldots,a_n|V) := f_{|V|}(a_{i_1},\ldots,a_{i_\ell})$ if $V = \{i_1<\cdots< i_\ell\}$.
\begin{defi}
The \textit{free} (\cite{Spe94}), \textit{Boolean} (\cite{SpW}), \textit{and monotone} (\cite{HS11}) \textit{functionals cumulants} form, respectively, the families of multilinear functionals $\{r_n:\A^n\to \mathbb{C}\}_{n\geq1}$, $\{b_n:\A^n\to \mathbb{C}\}_{n\geq1}$, and $\{h_n:\A^n\to \mathbb{C}\}_{n\geq1}$, implicitly defined by the equations
\begin{eqnarray}
\varphi(a_1\cdots a_n) &=& \sum_{\pi\in\NC(n)} r_{\pi}(a_1,\ldots,a_n),
\\\varphi(a_1\cdots a_n) &=& \sum_{\pi\in\operatorname{I}(n)}  b_{\pi}(a_1,\ldots,a_n),
\\\varphi(a_1\cdots a_n) &=& \sum_{\pi\in\NC(n)} \frac{1}{t(\pi)!} h_{\pi}(a_1,\ldots,a_n),
\end{eqnarray}
for any $n\geq1$ and $a_1,\ldots,a_n\in \A$. 
\end{defi}
A natural question arises when we want to find combinatorial formulas that relate the different brands of cumulants. The answer to this problem was initiated in \cite{Lehner} and studied in detail in \cite{AHLV}  using the standard approach based on Möbius inversion in several lattices of set partitions. 

\begin{theo}[\cite{AHLV}]
Let $\{r_n\}_{n\geq1}, \{b_n\}_{n\geq1}$, and $\{h_n\}_{n\geq1}$ be the be the families of free cumulants, Boolean cumulants, and monotone cumulants, respectively, in the non-commutative probability space $(\A,\varphi)$. Then, for any $n \geq 1$ and elements $a_1,\ldots,a_n \in \A$ we have that
\begin{eqnarray}
b_n(a_1,\ldots,a_n) &=& \sum_{\pi\in \NC_{\mathrm{irr}}(n)} r_\pi(a_1,\ldots,a_n),
\\ r_n(a_1,\ldots,a_n) &=& \sum_{\pi\in \NC_{\mathrm{irr}}(n)} (-1)^{|\pi|-1 }b_\pi(a_1,\ldots,a_n),\\
 b_n(a_1,\ldots,a_n) &=& \sum_{\pi\in \NC_{\mathrm{irr}}(n)} \frac{1}{t(\pi)!} h_\pi(a_1,\ldots,a_n), \label{eq:MonToBool}
\\ r_n(a_1,\ldots,a_n) &=& \sum_{\pi\in \NC_{\mathrm{irr}}(n)}  \frac{(-1)^{|\pi|-1}}{t(\pi)!} h_\pi(a_1,\ldots,a_n). \label{eq:MonToFree}
\end{eqnarray}
\end{theo}
The picture in \cite{AHLV}, however, was not complete since the formulas that give monotone cumulants in terms of free and Boolean cumulants were missing in the multivariate case. An answer to this issue was given in \cite{cefpp21}, where the authors found that the coefficients that describe the transition from free and Boolean cumulants to monotone cumulants are given by the Murua's coefficients of Definition \ref{def:murua}.

\begin{theo}[{\cite[Theorem 3]{cefpp21}}]
Let $\{r_n\}_{n\geq1}, \{b_n\}_{n\geq1}$, and $\{h_n\}_{n\geq1}$ be the be the families of free cumulants, Boolean cumulants, and monotone cumulants, respectively, in the non-commutative probability space $(\A,\varphi)$. Then, for any $n \geq 1$ and elements $a_1,\ldots,a_n \in \A$ we have that
\begin{eqnarray}
\label{eq:BoolToMon}
h_n(a_1,\ldots,a_n) &=& \sum_{\pi\in \NCirr(n)} \omega(t(\pi)) b_\pi(a_1,\ldots,a_n),\\
h_n(a_1,\ldots,a_n) &=& \sum_{\pi\in \NCirr(n)} (-1)^{|\pi|-1}\omega(t(\pi)) r_\pi(a_1,\ldots,a_n).\label{eq:FreeToMon}
\end{eqnarray}
\end{theo}
The main strategy in \cite{cefpp21} was to consider the free, Boolean, and monotone functionals as linear forms $\nu,\beta,\rho:L \to\mathbb{C}$, respectively, where $L:=T(\A)= \bigoplus_{n>0} \A^{\otimes n}$. More precisely, if $\{r_n\}_{n\geq1}$ is the family of free cumulants in $(\A,\varphi)$, we define the linear form $\nu(a_1\cdots a_n) := r_n(a_1,\ldots,a_n)$ for any word $a_1\cdots a_n\in \A^{\otimes n}$, and analogously for the Boolean and monotone cumulants. 
\par On the other hand, we know that $L$ is a pre-Lie algebra with respect to the pre-Lie product $\lhd$ described in \eqref{eq:pre-LieWords}. One of the main results of the series of papers \cite{EFP15, EFP16, EFP18, EFP19}, where the authors introduced a Hopf-algebraic framework for non-commutative probability, establishes that the relations monotone-free and monotone-Boolean cumulants are given through the pre-Lie exponential and the pre-Lie Magnus expansion in $L$.

\begin{theo}
\label{theo:MagnusRelations}
The linear forms associated to the monotone, free, and Boolean cumulants in $(\A,\varphi)$, denoted respectively by $\rho, \nu, \beta: T(\A)\to\mathbb{C}$, are related in terms of the pre-Lie exponential
\begin{eqnarray}
\label{eq:rhotoBool}
\beta &=& \exp^{\lhd}(\rho),\\ \label{eq:rhotokappa}\nu &=& -\exp^{\lhd}(-\rho),
\end{eqnarray}
and the pre-Lie Magnus expansion
\begin{equation}
\label{eq:kappatorho}
    \rho = \Omega(\beta) = -\Omega(-\nu).
\end{equation}
\end{theo}

\par Our goal in the rest of this section is to obtain the monotone-free and monotone-Boolean cumulant relations by applying the machinery developed earlier in the article. Given the pair $(\A,\varphi)$ consider the pre-Lie algebra of words over $\A$, $L = T(\A)$ equipped with the product \eqref{eq:pre-LieWords}. Recall that we can identify $L$ with its graded dual $L^*$. In this case, $L^*$ has a countable basis $\mathcal{B} = X^* = \{w_i\}_{i\in\mathbb{N}^*}$ given by the collection of non-empty words on the alphabet $\A$. Also, recall that the coproduct of $\kappa[L^*] $ given in Proposition \ref{prop:dualcop} can be written as
$$\overline{\delta}(w_i) = \sum_{i_0,I\neq \emptyset} \lambda^{i;i_0}_{I}w_{i_0}\otimes w_I$$
where if $I = \{i_1,i_2\ldots,i_s\}$, then $w_I = w_{i_1}\cdot w_{i_2}\cdot  \ldots\cdot w_{i_s}$, where the $\cdot$ stands for the product of polynomials in $\kappa[L^*]$. According to the definition of $\overline{\delta}$ in \eqref{eq:newdualcop}, one can easily obtain the following description of the $\lambda$ coefficients.
\begin{lemma}
\label{lem:auxS7}
Let $\{w_i\}_{i\in\mathbb{N}^*}$ be the set of non-empty words on $\A$. For any indexes $i,i_0,i_1,\ldots,i_s$, the coefficient $\lambda^{i;i_0}_{i_1,\ldots,i_s}$ is the number of ways in which it is possible to select subwords $v_1,v_2,\ldots,v_{s+1}$ of $w_i$ with $v_1,v_{s+1}\neq \emptyset$, such that $w_i = v_1w_{i_{\sigma(1)}}v_2 w_{i_{\sigma(2)}}\cdots v_sw_{i_{\sigma(s)}}v_{s+1}$ and $w_{i_0} = v_1v_2\cdots v_{s+1}$, where $\sigma$ is an arbitrary permutation in $S_s$.
\end{lemma}

With the purpose of applying the methods of the previous section to the pre-Lie algebra of words, it will be useful to relate the set of irreducible non-crossing partitions and the set of decorated trees that appear in the forest formula. This is done by the next lemma.

\begin{lemma}
\label{lem:NCForest}
Let $w_i$ be a non-empty word on $\A$ and $k\geq1$. Denote by $\T^{\prime,k}_{i}$ the subset of decorated trees associated to $w_i$, $T\in \T_i$ such that $\lambda(T)\neq 0$ and $|T|=k$. Then, for any $k\geq1$, there is a surjection $G:\NCirr^k(|w_i|)\to \T^{\prime,k}_i$ such that $|G^{-1}(T)| = \lambda(T)$, for any $T\in \T^{\prime,k}_i$.
\end{lemma}
\begin{proof}
First, we will give the description of the map $G$. For any $\pi\in 
\NCirr^k(|w_i|)$, we consider the tree of nestings of $\pi$, $t(\pi)$. It is cleat that $|t(\pi)| = k$. Assume that $w_i= a_1\cdots a_n$ with $a_1,\ldots,a_n\in\A$. Now, we make $t(\pi)$ a decorated tree $T$ as follows: for any $x\in V(t(\pi))$, if $x$ is associated to the block $V = \{j_1<j_2<\cdots < j_\ell\}\in\pi$, then $d_1(x)$ is the index in the list of words $\mathcal{B}=\{w_j\}_{j\in\mathbb{N}^*}$ such that $w_{d_1(x)} = a_{j_1}a_{j_1+1}\cdots a_{j_2}a_{j_2+1}\cdots a_{j_\ell}$, and $d_2(x)$ is the index such that $w_{d_2(x)} = a_{j_1}a_{j_2}\cdots a_{j_\ell}$. Since clearly $d_1(\operatorname{root}(T)) = w_i$ and $\lambda^{d_1(x);d_2(x)}_{d_1(\operatorname{succ}(x))}\neq 0$ for any $x\in \operatorname{Int}(T)$, then $T\in \T^{\prime,k}_i$. We can define then $G:\NCirr^k(|w_i|)\to \T^{\prime,k}_i$ by $G(\pi) = T$, for any $\pi\in \NCirr^k(|w_i|)$.
\par Now, we will prove that the map $G$ is surjective. Take $T\in \T^{\prime,k}_i$. If $T$ is a single-vertex tree, i.e. $k=1$, we consider the single-block partition $1_{|w_i|} \in \NCirr^1(|w_i|)$.  Otherwise, assume that $T = B^+_{(i,i_0)}(T_1,\ldots,T_s)$, where $T_j$ is associated to $w_{i_j}$ for any $1\leq j\leq s$. Since $\lambda(T)\neq 0$, $\lambda^{i,i_0}_{i_1,\ldots,i_s}= m>0$. By Lemma \ref{lem:auxS7}, there are $m$ ways to take $s+1$ subwords $v_1,\ldots,v_{s+1}$ of $w_i = a_1\cdots a_n$ such that $v_1,v_{s+1}\neq\emptyset$, $w_{i_0}= v_1\cdots v_{s+1}$, and $w_i = v_1w_{i_{\sigma(1)}}v_2w_{i_{\sigma(2)}}\cdots v_sw_{i_{\sigma(s)}}v_{s+1}$ for a permutation $\sigma\in S_s$. Each of the $m$ possible selection of the subwords $v_1,\ldots,v_{s+1}$ corresponds to a selection of pairwise disjoint subsets $C_1,\ldots,C_{s+1}\subseteq [n]$ such that $1\in C_1$ and $n\in C_{s+1}$. We now construct the block $V_0 = C_1\cup \cdots \cup C_{s+1}$. The remaining indexes $[n]\backslash V_0$ can be grouped into $s$ non-empty pairwise disjoint intervals $J_1,\ldots,J_s$ such that $J_l = \{r_l,r_l+1,\ldots,r_l+m\}$ is the subset of indexes such that $w_{i_{\sigma(l)}} = a_{r_l}a_{r_l+1}\cdots a_{r_l+m}$, for any $1\leq l \leq s$.
Due to the possible repetitions in the subwords $w_{i_j}$, there are $\operatorname{sym}(B^-(T))$ ways to allocate the decorated trees $T_1,\ldots,T_s$ to the subsets $J_1,\ldots,J_s\in [n]\backslash V_0$. 
\par Finally, we proceed by induction. Indeed, since for any $1\leq l\leq s$, $T_l$ is associated to a $J_p$ and $|T_l| < k$, we can find $\lambda(T_l)$ different irreducible non-crossing partitions $\pi_l \in \NCirr^{|T_l|}(J_p)$  such that their corresponding decorated tree of nestings is $T_l$. Therefore, we can construct $\lambda^{i;i_0}_{i_1,\ldots,i_s}\operatorname{sym}(B^-(T))\lambda(T_1)\cdots \lambda(T_s) = \lambda(T)$ irreducible non-crossing partitions of the form
$$\pi = \{V_0\}\cup \pi_1\cup\cdots\cup\pi_s \in \NCirr(n)$$
such that $1,n\in V_0$, $|\pi| = 1 + |T_1| + \cdots + |T_s| = |T| = k$ and $G(\pi) = T$. This completes the proof. 
\end{proof}

The following two theorems were obtained in \cite{cefpp21} by direct computations, we show here how they can be deduced from pre-Lie forest formulas.

\begin{theo}
\label{thm:exp}
Let $L$ be the pre-Lie algebra of words associated to an alphabet $X$. Then, for $\alpha\in L$ and a word $w_i \in X^*\subset L^*$ such that $w_i= a_1\cdots a_n$ with $a_1,\ldots,a_n\in X$, we have that
$$\langle \exp^{\lhd}(\alpha) | w_i\rangle = \sum_{\pi\in \NC_{\mathrm{irr}}(n)} \frac{1}{t(\pi)!} \alpha_\pi(w_i),$$
where $\alpha_\pi(w_i) =  \prod\limits_{V\in \pi} \langle \alpha | w_V\rangle$, where $w_V = a_{j_1}a_{j_2}\cdots a_{j_r}$ if $V = \{j_1< j_2 < \cdots < j_r\}$.
\end{theo}

\begin{proof}
By the definition of the pre-Lie exponential, Lemma \ref{lem:iterated2} and the forest formula \eqref{eq:forestformula}, we have that $  
\langle\exp^{\lhd}(\alpha)|w_i\rangle = \sum_{k\geq1} \frac{1}{k!} \langle r^{(k)}_{\alpha}(\alpha) |w_i\rangle$ and

\begin{eqnarray*}
\langle  r^{(k)}_{\alpha}(\alpha) |w_i\rangle &=& \langle \alpha\otimes \cdots \otimes \alpha | \delta_{\mathrm{irr}}^{[k]} (w_i)\rangle
\\ &=& \left\langle \alpha\otimes \cdots \otimes \alpha \left| \sum_{T\in \mathcal{T}_i} \sum_{f\in \mathrm{lin}(T)} \lambda(T) C(f) \right. \right\rangle
\\ &=& \sum_{\substack{T\in \mathcal{T}_i^\prime \\ |T| = k}}\lambda(T) \sum_{\substack{f\in \mathrm{lin}(T)}}  \prod_{j=1}^k  \langle \alpha | w_{d_2(f^{-1}(j))}\rangle 
\\ &=& \sum_{\substack{T\in \mathcal{T}_i^\prime \\ |T| = k}} \lambda(T)m(T) \prod_{x\in V(T)} \langle \alpha | w_{d_2(x)}\rangle,
\end{eqnarray*}
where $\T'_i$ stands for the subset of decorated trees associated to $w_i$, $T\in \T_i$, such that $\lambda(T)\neq 0$. Thus, we can apply Lemma \ref{lem:NCForest} in order to rearrange the sum in the right-hand side of the last equation as:

\begin{equation}
\label{eq:proof73}
\langle \alpha\otimes\cdots\otimes \alpha | \delta_{\mathrm{irr}}^{[k]}(w_i)\rangle = \sum_{\pi\in \NCirr^k(n)} m(t(\pi)) \prod_{V\in \pi}\langle \alpha | w_V\rangle.
\end{equation}

The above equation says, in particular, that $k$ cannot be greater or equal than $n$ since an irreducible non-crossing partition $\pi\in \NCirr(n)$ has at most $n-1$ blocks. Hence the pre-Lie exponential can be written in the following way:

\begin{eqnarray*}
\langle \exp^\lhd(\alpha) | w_i \rangle &=& \sum_{k=1}^\infty \frac{1}{k!} \langle \alpha\otimes\cdots \otimes \alpha | \delta_{\mathrm{irr}}^{[k]}(w_i)\rangle
\\ &=& \sum_{k=1}^{n-1} \frac{1}{k!} \sum_{\pi\in \NC_{\mathrm{irr}}^k(n)} m(t(\pi)) \alpha_\pi(w_i)
\\ &=&  \sum_{k=1}^{n-1}\sum_{\pi\in \NC_{\mathrm{irr}}^k(n)} \frac{m(t(\pi))}{|\pi|!} \alpha_\pi(w_i)
\\ &=& \sum_{\pi\in \NC_{\mathrm{irr}}(n)} \frac{1}{t(\pi)!} \alpha_\pi(w_i),
\end{eqnarray*}
as we wanted to show.
\end{proof}

\begin{rem}
In the context of the above theorem, by considering $\alpha = \rho$ the linear form associated to the monotone cumulants in $(\A,\varphi)$, by the previous theorem and \eqref{eq:rhotoBool}, we obtain the monotone-to-Boolean cumulant relation \eqref{eq:MonToBool}. On the other hand, by considering $\alpha = -\rho$, from \eqref{eq:rhotokappa} we obtain the monotone-to-free cumulant relation \eqref{eq:MonToFree}.
\end{rem}

The final theorem establishes a combinatorial formula for the Magnus' expansion in the pre-Lie algebra of words. In particular, it allows us to recover the free-to-monotone and Boolean-to-monotone cumulant formulas for the multivariate case proved in \cite{cefpp21}.

\begin{theo}
\label{thm:MagnusWords}
Let $L$ be the pre-Lie algebra of words associated to an alphabet $X$. Then, for $\alpha\in L$ and a word $w_i\in X^*\subset L^*$ such that $w_i= a_1\cdots a_n$ with $a_1,\ldots,a_n\in X$, we have that
\begin{eqnarray}
\label{eq:MagnusWords}
\langle \Omega(\alpha)|w_i \rangle = \sum_{\pi\in \NCirr(n)} \omega(t(\pi)) \alpha_\pi(w_i).
\end{eqnarray}
\end{theo}

\begin{proof}
We will follow the strategy used in the proof of Proposition \ref{prop:MagnusFree} in Section \ref{sec:MuruaMagnusI}. By Theorem \ref{thm:AltMagnus}, we have the following computation for $\alpha\in L$ and $w_i=a_1\cdots a_n \in L^*$:
\begin{eqnarray*}
\langle \Omega(\alpha) | w_i\rangle &=& \langle \sol(\exp(\alpha)) | w_i \rangle
\\ &=& \sum_{m\geq0} \frac{1}{m!}\langle \sol(\alpha^{\cdot m}) | w_i \rangle
\\ &=& \sum_{m\geq0} \frac{1}{m!}\sum_{k=1}^m \frac{(-1)^{k-1}}{k} \sum_{\substack{j_1,j_2,\ldots,j_k\geq1\\ j_1+ j_2 + \cdots + j_k = m}} {m\choose j_1,j_2,\ldots,j_k} \langle \alpha^{j_1}* \alpha^{j_2}* \cdots * \alpha^{j_k} | w_i\rangle
\end{eqnarray*}
By the same arguments that in the proof of Proposition \ref{prop:MagnusFree}, we have that only the case $j=1$ may produce a non-zero contribution in the right-hand side of the above equation. Also, for the same reason, the iterated $*$ products can be replaced by iterated brace products (see \cite[Chap. 6]{CP21} for details on such arguments), and by Lemma \ref{lem:iteratedbrace} we obtain
$$ \langle \Omega(\alpha) | w_i\rangle = \sum_{m\geq 0 }\sum_{k=1}^m \frac{(-1)^{k-1}}{k} \sum_{\substack{j_2,\ldots,j_k\geq1\\1+j_2+\cdots+j_k = m}} \frac{1}{j_2!\cdots j_k!}\langle \alpha \otimes \alpha^{j_2}\otimes\cdots \otimes \alpha^{j_k} \, | \, \overline\delta^{[k]}(w_i)\rangle.$$
Observe that, by the definition of the reduced coproduct, a $m>n$ will produce a zero contribution in $\langle \Omega(\alpha)|w_i\rangle$. Hence, by the forest formula for the iterated reduced coproduct \eqref{eq:forestformula0}, we have
\begin{eqnarray*}
\langle \Omega(\alpha)|w_i\rangle &=& \sum_{m=1}^n \sum_{k=1}^m \frac{(-1)^{k-1}}{k} \sum_{T\in \T'_i} \sum_{f\in k-\operatorname{lin}(T)} \lambda(T) \sum_{\substack{j_2,\ldots,j_k\geq1\\1+j_2+\cdots+j_k = m}} \langle \alpha \otimes \alpha^{j_2}\otimes\cdots \otimes \alpha^{j_k} \, | \, C(f)\rangle,
\end{eqnarray*}
where $\T_i'$ stands for the subset of decorated trees $T\in \T_i$ associated to $w_i$ such that $\lambda(T)\neq 0$. Now, given $1\leq m\leq n$, $1\leq k\leq m$, and a tuple $(j_2,\ldots,j_k)$ such that $1+j_2+\cdots + j_k=m$, the only decorated trees $T\in \T '_i$ such that there exists a $f\in k-\operatorname{lin}(T)$ which produces a non-zero contribution for $ \langle \alpha \otimes \alpha^{j_2}\otimes\cdots \otimes \alpha^{j_k} \, | \, C(f)\rangle$ must satisfy that $|T| = m$. On the other hand, for $T\in \T_i'$ and $f\in k-\operatorname{lin}(T)$ we have
\begin{eqnarray*}
\langle \alpha \otimes \alpha^{j_2}\otimes\cdots \otimes \alpha^{j_k} | C(f) \rangle &=& \langle \alpha| w_{d_2(f^{-1}(1))}\rangle \prod_{l=2}^k \langle \alpha^{j_l} | w_{d_2(f^{-1}(l))} \rangle
\\ &=& \langle \alpha| w_{d_2(f^{-1}(1))}\rangle \prod_{l=2}^k \left( \sum_{\sigma \in S_{j_l}} \prod_{h\in f^{-1}(l)} \langle \alpha | w_h \rangle\right)
\\ &=& j_2!\cdots j_k! \prod_{x\in V(T)} \langle \alpha| w_{d_2(x)}\rangle,
\end{eqnarray*}
where in the second equality, we used the duality described in \eqref{eq:extension}. Combining this with the fact that, given $f\in k-\operatorname{lin}(T)$, there is exactly one tuple $(j_2',\ldots,j_k')$ of positive integers such that $\langle \alpha \otimes \alpha^{j_2'}\otimes\cdots \otimes \alpha^{j_k'} | C(f) \rangle\neq 0$, and that it is given by $j_i'= |f^{-1}(l)|$ for $2\leq l\leq k$, it follows that
\begin{eqnarray*}
\langle \Omega(\alpha) |w_i \rangle &=& \sum_{m=1}^n \sum_{k=1}^m \frac{(-1)^{k-1}}{k} \sum_{\substack{T\in \T'_i\\ |T|=m}} \lambda(T) \sum_{f\in k-\operatorname{lin}(T)}  \prod_{x\in V(T)} \langle \alpha| w_{d_2(x)}\rangle
\\ &=&  \sum_{m=1}^n \sum_{k=1}^m \frac{(-1)^{k-1}}{k} \sum_{\substack{T\in \T'_i\\ |T|=m}} \lambda(T)  \prod_{x\in V(T)} \langle \alpha| w_{d_2(x)}\rangle \sum_{f\in k-\operatorname{lin}(T)}  1
\\ &=& \sum_{m=1}^n \sum_{k=1}^m \frac{(-1)^{k-1}}{k} \sum_{\substack{T\in \T'_i\\ |T|=m}} \lambda(T)  \prod_{x\in V(T)} \langle \alpha| w_{d_2(x)}\rangle \omega_k(T),
\end{eqnarray*}
where we recall that $\omega_k(T) = |k-\operatorname{lin}(T)|$. As in the proof of Theorem \ref{thm:exp}, we use Lemma \ref{lem:NCForest} in order to rewrite the above equation as
\begin{eqnarray*}
\langle \Omega(\alpha)|w_i\rangle &=& \sum_{m=1}^n \sum_{k=1}^m \frac{(-1)^{k-1}}{k} \sum_{\pi\in \NCirr^m(n)}  \omega_k(t(\pi)) \prod_{V\in \pi} \langle \alpha| w_{V}\rangle
\\ &=& \sum_{m=1}^n \sum_{\pi\in \NCirr^m(n)} \alpha_\pi(w_i) \sum_{k=1}^m \frac{(-1)^{k-1}}{k}\omega_k(t(\pi))
\\ &=& \sum_{\pi\in \NCirr(n)} \omega(t(\pi)) \alpha_\pi(w_i),
\end{eqnarray*}
which concludes the proof.
\end{proof}

\begin{rem}
In the context of Theorem \ref{thm:MagnusWords}, by considering $\alpha = \beta$ the linear form associated to the Boolean cumulants in $(\A,\varphi)$, by Theorem \ref{thm:MagnusWords} and \eqref{eq:kappatorho}, we obtain the Boolean-to-monotone cumulant relation \eqref{eq:BoolToMon}. On the other hand, by considering $\alpha = -\nu$, where $\nu$ is the linear form associated to the free cumulants,  from \eqref{eq:kappatorho} we obtain the free-to-monotone cumulant relation \eqref{eq:FreeToMon}.
\end{rem}

As we did in Section \ref{sec:MuruaMagnusII} for the action of the Magnus operator in the free pre-Lie algebra, we can give an alternative proof of Theorem \ref{thm:MagnusWords} by using the recursion \eqref{eq:MuruaRecursion} for the Murua's coefficients. We include it for completeness sake, as it provides further insights on the combinatorics of free probability.

\begin{proof}[Alternative proof of Theorem \ref{thm:MagnusWords}]
Once more, by definition of the pre-Lie  Magnus operator \eqref{eq:MagnusExpansion}, Lemma \ref{lem:iterated2}, and the forest formula \eqref{eq:forestformula}, we have that 
\begin{eqnarray*}
\langle \Omega(\alpha) | w_i\rangle &=& \sum_{m\geq1} \frac{B_m}{m!}\langle r^{(m)}_{\Omega(\alpha)}(\alpha)|w_i \rangle
\\ &=& \sum_{m\geq1} \frac{B_m}{m!}\langle \alpha \otimes \Omega(\alpha)\otimes\cdots \otimes \Omega(\alpha) | \delta_{\mathrm{irr}}^{[m+1]}(w_i)\rangle
\\ &=& \sum_{m=1}^{n-1} \frac{B_m}{m!} \sum_{\substack{T\in \T'_i\\ |T| = m+1}} \sum_{f\in\operatorname{lin}(T)} \lambda(T)\langle \alpha \otimes \Omega(\alpha)\otimes\cdots \otimes \Omega(\alpha) | C(f) \rangle,
\end{eqnarray*}
where we recall that $\T'_i$ is the subset of decorated trees associated to $w_i$, $T\in \T_i$, such that $\lambda(T)\neq 0$. By the same argument that in the proof of Proposition \ref{prop:MagnusFree}, the term $\langle \alpha \otimes \cdots \otimes \Omega(\alpha)| C(f)\rangle$ does not depend of $f$. Hence, by using Lemma \ref{lem:NCForest}, we get that
\begin{eqnarray*}
\langle\Omega(\alpha) |w_i \rangle &=& \sum_{T\in \T'_i} \frac{B_{|T|-1}}{(|T|-1)!} m(T) \lambda(T) \langle \alpha | w_{d_2(\operatorname{root}(T))}\rangle \prod_{x\in V(T)\backslash \{\operatorname{root}(T)\}} \langle \Omega(\alpha) | w_{d_2(x)}\rangle
\\ &=& \sum_{\substack{\pi \in \NCirr(n)\\ 1,n\in V_0}} \frac{B_{|\pi|-1}}{t(\pi\backslash \{V_0\})!} \langle \alpha | w_{V_0} \rangle \prod_{\substack{W\in \pi\\W\neq V_0}} \langle \Omega(\alpha)|w_W\rangle.  
\end{eqnarray*}
Since $|w_W|< n$ for any $W\in \pi$, we can conclude by induction in the same way that in the second part of the proof of Theorem 3 in \cite{cefpp21}. We will outline the ideas for the convenience of the reader. Using induction on each $\langle \Omega(\alpha)|w_W\rangle$, we get
$$ \langle \Omega(\alpha)|w_i\rangle = \sum_{\substack{\pi \in \NCirr(n)\\ 1,n\in V_0}} \frac{B_{|\pi|-1}}{t(\pi\backslash \{V_0\})!} \langle \alpha | w_{V_0} \rangle \prod_{\substack{W\in \pi\\W\neq V_0}} \left(\sum_{\sigma_W \in\NCirr(W)} \omega(t(\sigma_W)) \alpha_{\sigma_W}(w_W) \right).$$
Notice that $\{V_0\}\cup\bigsqcup_{W\in \pi\backslash\{V_0\}} \sigma_W$ is a non-crossing partition. 
There is a bijection between 
\begin{itemize}
    \item pairs $(\pi\backslash\{V_0\},(\sigma_W)_{W\in \pi\backslash\{V_0\}})$, where $\pi$ is a non-crossing irreducible partition, $V_0$ is its outer block (the block which contains $1$ and $n$) and the $\sigma_W$ are  non-crossing irreducible partitions of its other blocks, and
    \item pairs $(\mu',S)$, where $\mu=\mu'\cup V_0$ is an irreducible non-crossing partition with outer block $V_0$ and $S$ is a subset of the set of blocks of $\mu'$.
\end{itemize}
One sets $\mu' := \bigsqcup_{W\in \pi\backslash\{V_0\}} \sigma_W$ and, if $V_{0,W}$ stands for the unique outer block of $\sigma_W$, one considers the collection $S:= \{V_{0,W}\}_{W\in \pi\backslash \{V_0\}}$. Conversely, by an argument that we omit, one can reconstruct from $(\mu',S)$ the blocks $W$ (they are in bijection with $S$) and the $\sigma_W$ (they are subsets of the set of blocks of $\mu'$) using the poset structure on the blocks of $\mu'$ and the induced poset structure on $S$.
In the language of the notation of Proposition \ref{prop:recursiveMurua}, it is possible to show that $S\in K(t(\mu'))$ and $\{t(\sigma_W)\}_{W\in \pi\backslash \{V_0\}} \in C^S(t(\mu'))$. Thus, using the previous bijection and $\mu = \mu'\cup \{V_0\}$, we obtain that
\begin{eqnarray*}
\langle \Omega(\alpha)|w_i\rangle &=& \sum_{\substack{\mu\in \NCirr(n)\\1,n\in V_0}} \sum_{S\in K(t(\mu'))} \frac{B_{|S|}}{S!} \langle \alpha| w_{V_0} \rangle \left( \prod_{W\in \mu'} \langle \alpha | w_W\rangle \right) \omega(C^S(t(\mu')))
\\ &=& \sum_{\mu\in \NCirr(n)} \omega(t(\mu))\alpha_\mu(w_i),
\end{eqnarray*}
where in the last equality we used Proposition \ref{prop:recursiveMurua}. This concludes the proof.
\end{proof}

\begin{rem}
The proof of Theorem 3  in \cite{cefpp21} uses fundamentally Proposition 4 therein, which gives a combinatorial expression in terms of monotone non-crossing partitions for the iterated pre-Lie product. In the approach explained in this work, the role of Proposition 4 is taken by the forest formula \eqref{eq:forestformula}. 
\end{rem}

\subsection*{Acknowledgements} The first author was partially supported
by the project Pure Mathematics in Norway, funded by Trond Mohn Foundation and Tromsø Research Foundation. The second author acknowledges support from the European Research Council (ERC) under the European Union’s
Horizon 2020 research and innovation program (Duall project, grant agreement No. 670624) and from the ANR project Algebraic Combinatorics, Renormalization, Free probability and Operads – CARPLO (Project-ANR-20-CE40-0007).



\end{document}